\documentclass[12pt]{amsart}
\usepackage{amsmath,amssymb,latexsym,esint,cite,mathrsfs}
\usepackage{verbatim,wasysym}
\usepackage[left=2.6cm,right=2.6cm,top=3cm,bottom=3cm]{geometry}
\makeatletter \@addtoreset{equation}{section} \makeatother

\usepackage{graphicx}
\usepackage{subfigure}
\usepackage{graphicx,epic,eepic}
\usepackage{tikz}  %ͼÐλæÖƺê°ü
\usepackage[colorlinks=true,urlcolor=blue,
citecolor=red,linkcolor=blue,linktocpage,pdfpagelabels,
bookmarksnumbered,bookmarksopen]{hyperref}
\numberwithin{equation}{section}
\newtheorem{theorem}{Theorem}[section]
\newtheorem{lemma}[theorem]{Lemma}

\newtheorem{corollary}[theorem]{Corollary}
\newtheorem{remark}[theorem]{Remark}
\numberwithin{equation}{section}

\allowdisplaybreaks

\begin{document}

\title[Symmetry breaking of extremals for CKN inequality]
{Symmetry breaking of extremals for the high order Caffarelli-Kohn-Nirenberg type inequalities}

\author[S. Deng]{Shengbing Deng}
\address{\noindent Shengbing Deng
\newline
School of Mathematics and Statistics, Southwest University,
Chongqing 400715, People's Republic of China}\email{shbdeng@swu.edu.cn}

\author[X. Tian]{Xingliang Tian$^{\ast}$}
\address{\noindent Xingliang Tian  \newline
School of Mathematics and Statistics, Southwest University,
Chongqing 400715, People's Republic of China.}\email{xltian@email.swu.edu.cn}

\thanks{$^{\ast}$ Corresponding author}

\thanks{2020 {\em{Mathematics Subject Classification.}} 35J30, 49K30, 26D10}

\thanks{{\em{Key words and phrases.}} Caffarelli-Kohn-Nirenberg inequalities; Weighted fourth-order equation; Non-degeneracy; Symmetry and symmetry breaking}

\allowdisplaybreaks
%\maketitle

\begin{abstract}
{\tiny
In this paper we give the first result about the precise symmetry and symmetry breaking regions of extremal functions for weighted second-order inequalities. Firstly, based on the work of C.-S. Lin [Comm. Partial Differential Equations, 1986], a new second-order Caffarelli-Kohn-Nirenberg type inequality will be established, i.e., 
    \begin{equation*}
    \int_{\mathbb{R}^N}|x|^{-\beta}|\mathrm{div} (|x|^{\alpha}\nabla u)|^2 \mathrm{d}x
    \geq \mathcal{S}\left(\int_{\mathbb{R}^N}
    |x|^{\beta}|u|^{p^*_{\alpha,\beta}} \mathrm{d}x\right)^{\frac{2}{p^*_{\alpha,\beta}}},\quad \mbox{for all}\ u\in C^\infty_0(\mathbb{R}^N),
    \end{equation*}
    for some constant $\mathcal{S}=\mathcal{S}(N,\alpha,\beta)>0$, where
    \begin{align*}
    N\geq 5,\quad \alpha>2-N,\quad \alpha-2<\beta\leq \frac{N}{N-2}\alpha,\quad p^*_{\alpha,\beta}=\frac{2(N+\beta)}{N-4+2\alpha-\beta}.
    \end{align*}
    We obtain a symmetry breaking conclusion: when $\alpha>0$ and $\beta_{\mathrm{FS}}(\alpha)<\beta< \frac{N}{N-2}\alpha$ where $\beta_{\mathrm{FS}}(\alpha):=
        -N+\sqrt{N^2+\alpha^2+2(N-2)\alpha}$, then the extremal function for the best constant $\mathcal{S}$, if it exists, is nonradial. Furthermore, we give a symmetry result when $\beta=\frac{N}{N-2}\alpha$ and $2-N<\alpha<0$, i.e., for all $u\in C^\infty_0(\mathbb{R}^N)$,
        \begin{equation*}
    \int_{\mathbb{R}^N}|x|^{-\frac{N}{N-2}\alpha}|\mathrm{div} (|x|^{\alpha}\nabla u)|^2 \mathrm{d}x
    \geq \left(1+\frac{\alpha}{N-2}\right)^{4-\frac{4}{N}}\mathcal{S}_0
    \left(\int_{\mathbb{R}^N}
    |x|^{\frac{N}{N-2}\alpha}|u|^{\frac{2N}{N-4}} \mathrm{d}x\right)^{\frac{N-4}{N}},
    \end{equation*}
    where $\mathcal{S}_0$ is the sharp constant of second-order Sobolev inequality, and equality holds only for some radial functions of the form $(1+|x|^{2+\frac{2\alpha}{N-2}})
    ^{-\frac{N-4}{2}}$ (up to scalings and multiplications).
    }
\end{abstract}

\vspace{3mm}

\maketitle

\section{{\bfseries Introduction}}\label{sectir}

\subsection{Motivation}\label{subsectmot}
Let us recall the famous so-called first order  Caffarelli-Kohn-Nirenberg (we write (CKN) for short) inequality which was first introduced in 1984 by Caffarelli, Kohn and Nirenberg in their celebrated work \cite{CKN84}, here, we only focus on the case without interpolation term, that is,
    \begin{equation}\label{cknwit}
    \left(\int_{\mathbb{R}^N}|x|^{-b\tau}|u|^\tau \mathrm{d}x\right)^{\frac{p}{\tau}}
    \leq C_{\mathrm{CKN}}\int_{\mathbb{R}^N}|x|^{-ap}|\nabla u|^p \mathrm{d}x, \quad \mbox{for all}\ u\in C^\infty_0(\mathbb{R}^N),
    \end{equation}
    for some constant $C_{\mathrm{CKN}}>0$, where
    \begin{equation}\label{cknwitc}
    1<p<N,\quad -\infty<a<a_c:=\frac{N-p}{p},\quad a\leq b\leq a+1,\quad \tau=\frac{pN}{N-p(1+a-b)}.
    \end{equation}
    A natural question is whether the best constant $C_{\mathrm{CKN}}$ could be achieved or not? Moreover, if $C_{\mathrm{CKN}}$ is achieved, are the extremal functions radial symmetry? In fact, the weights $|x|^{-ap}$ and $|x|^{-b\tau}$ have deep influences in many aspects about this inequality, for example, achievable of best constant and symmetry of minimizers.

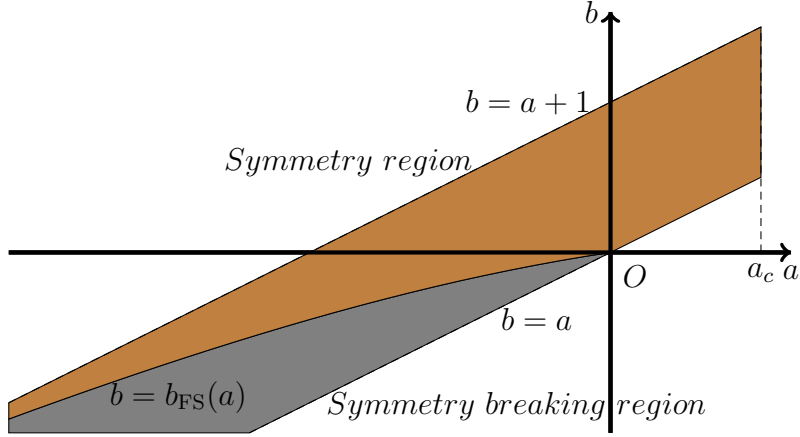
\begin{figure}
\begin{tikzpicture}[scale=4]
%%================ÉèÖÃ×ø±êÖ᷶Χ
		\draw[->,ultra thick](-2,0)--(0,0)node[below right]{$O$}--(0.6,0)node[below]{$a$};
		\draw[->,ultra thick](0,-0.6)--(0,0.8)node[left]{$b$};
%%===============»­Í¼
        \draw[fill=gray,domain=0:-2]plot(\x,{4*(1-0.47*\x)/
        (2*((1-0.6*\x)^2+3)^0.5)
        +0.47*\x -1})--(-2,-0.6)--(-1.2,-0.6)--(0,0);
        \draw[fill=brown,domain=0:-2]plot(\x,{4*(1-0.47*\x)/
        (2*((1-0.6*\x)^2+3)^0.5)
        +0.47*\x -1})--(-2,-0.5)--(0.5,0.75)--(0.5,0.25)--(0,0);
        \draw[densely dashed](0.5,0)node[below]{$a_c$}--(0.5,0.75);
        \draw[densely dashed](-2,-0.5)--(0.5,0.75);
        \draw[densely dashed](-1.2,-0.6)--(0,0);
        \draw[-,ultra thick](0,-0.6)--(0,0.8);
        \draw[-,ultra thick](-2,0)--(0.6,0);
%%===============ÉèÖÃͼÀý========================================================

		\node[left] at(-0.03,0.5){$b=a+1$};
        \node[right] at (-0.4,-0.22){$b=a$};
        \node[right] at (-1.7,-0.48){$b=b_{\mathrm{FS}}(a)$};
        \node[right] at (-0.98,-0.51){$Symmetry\ breaking\ region$};
        \node[right] at (-1.32,0.3){$Symmetry\ region$};
\end{tikzpicture}
\caption{\small For the first order case with $p=2$. The {\em Felli-Schneider region}, or symmetry breaking region, appears in dark grey and is defined by $\{(a,b): a<0,\ a<b<b_{\mathrm{FS}}(a)\}$. And symmetry holds in the brown region defined by $\{(a,b): a<0, \ b_{\mathrm{FS}}(a)\leq b<a+1\}\cup \{(a,b): 0\leq a<a_c, \ a\leq b<a+1\}$.}
\label{F1}
\end{figure}

    For $p=2$, there are complete results for the above problem. When $b=a+1$ or $a<0$ and $b=a$, Catrina and Wang \cite{CW01} proved that $C_{\mathrm{CKN}}$ is not achieved and for other cases it is always achieved. Furthermore, when $a<0$ and $a< b<b_{\mathrm{FS}}(a)$, where
    \begin{align}\label{deffsc}
    b_{\mathrm{FS}}(a):=\frac{N(a_c-a)}{2\sqrt{(a_c-a)^2+N-1}}
    +a-a_c,
    \end{align}
    Felli and Schneider \cite{FS03} proved the extremal function is non-radial by restricting it in the radial space and classifying linearized problem, thus $b_{\mathrm{FS}}(a)$ is usually called {\em Felli-Schneider curve}. Finally, in a celebrated paper, Dolbeault, Esteban and Loss \cite{DEL16} proved an optimal rigidity result by using the so-called {\em carr\'{e} du champ} method that when $a<0$ and $b_{\mathrm{FS}}(a)\leq b<a+1$, the extremal function is symmetry. %, and we refer to \cite{Do21} for an overall review about this method.
    See the previous results shown as in Figure \ref{F1}. For the general case $1<p<N$, we refer to \cite{BW02,CM13,LL17,CC22}, however, the optimal symmetry region is not known yet and it is conjectured that
    $\frac{(1+a-b)(a_c-a)}{(a_c-a+b)}
    \leq \sqrt{\frac{N-1}{N/(1+a-b)-1}}$.
    We also refer to \cite{ADG22,DGK24} for the fractional-order (CKN) inequality.

    In 1986, C.-S. Lin \cite{Li86} extended the (CKN) inequality of \cite{CKN84} to the higher-order case, here we only mention the second-order case without interpolation term:
    \vskip0.25cm

    \noindent{\bf Theorem~A.} \cite{Li86} {\it  There exists a constant $C>0$ such that
    \begin{equation}\label{ckn2y}
    \int_{\mathbb{R}^N}|x|^{-2a}|D^2 u|^2 \mathrm{d}x
    \geq C\left(\int_{\mathbb{R}^N}|x|^{-b\tau}|u|^\tau \mathrm{d}x\right)^{\frac{2}{\tau}},
    \quad \mbox{for all}\ u\in C^\infty_0(\mathbb{R}^N),
    \end{equation}
    where $-\infty<a<\frac{N-4}{2}$, $a\leq b\leq a+2$, $\tau=\frac{2N}{N-2(2+a-b)}$. Here,
    \[
    |D^2 u|^2=\sum\left|\frac{\partial^2 u}{\partial x_1^{m_1}\cdots \partial x_N^{m_N}}\right|^2,\ m_i\in \{0,1,2\} \ \mbox{satisfying}\ \sum_{i=1}^N m_i=2.
    \]
    }
    \vskip0.25cm

    Therefore, same as the first-order case, it is natural to ask that whether the best constant could be achieved or not? Moreover, if it is achieved, are the extremal functions radial symmetry? However, the second-order case \eqref{ckn2y} is complicated. After integration by parts we get $\||x|^{-a}D^2 u\|_{L^2(\mathbb{R}^N)}=\||x|^{-a}\Delta u\|_{L^2(\mathbb{R}^N)}$ unless $a=0$. When $N\geq 5$ and $0\leq a<\frac{N-4}{2}$, Szulkin and Waliullah \cite[Lemma 3.1]{SW12} proved that $\||x|^{-a}D^2 u\|_{L^2(\mathbb{R}^N)}$ is equivalent to $\||x|^{-a}\Delta u\|_{L^2(\mathbb{R}^N)}$. For general $a$, Caldiroli and Musina \cite[Remark 2.3]{CM11} proved that $\||x|^{-a}D^2 u\|_{L^2(\mathbb{R}^N)}$ is equivalent to $\||x|^{-a}\Delta u\|_{L^2(\mathbb{R}^N)}$ if and only if
 \begin{align}\label{defac}
 a\neq -\frac{N}{2}-k\quad\mbox{and}\quad a\neq \frac{N-4}{2}+k, \quad\mbox{for all}\ k\in\mathbb{N}.
 \end{align}
 Therefore, with addition condition \eqref{defac}, we have that \eqref{ckn2y} is equivalent to the following second-order CKN inequality
    \begin{equation}\label{ckn2}
    \int_{\mathbb{R}^N}|x|^{-2a}|\Delta u|^2 \mathrm{d}x
    \geq C_{\mathrm{CKN2}}\left(\int_{\mathbb{R}^N}|x|^{-b\tau}|u|^\tau \mathrm{d}x\right)^{\frac{2}{\tau}},
    \quad \mbox{for all}\ u\in C^\infty_0(\mathbb{R}^N),
    \end{equation}
    which has it's own interests.
    Furthermore, Caldiroli and Musina \cite[Remark 2.3]{CM11} also showed that \eqref{ckn2} does not hold if $a=-\frac{N}{2}-k$ or $\frac{N-4}{2}+k$ for some $k\in\mathbb{N}$, which comes from spectral analysis about the weighted Rellich inequality (see \cite{CM12}): for all $u\in C^\infty_0(\mathbb{R}^N)$,
    \begin{align}\label{wri}
    \int_{\mathbb{R}^N}|x|^{-2a}|\Delta u|^2 \mathrm{d}x
    \geq \inf_{k\in\mathbb{N}}
    \left(k+\frac{N}{2}+a\right)^2\left(k+\frac{N-4}{2}-a\right)^2
    \int_{\mathbb{R}^N}|x|^{-2a-4}|u|^2 \mathrm{d}x.
    \end{align}
    Now, let us focus on \eqref{ckn2} under condition \eqref{defac}. There are only some partial results about the existence of extremal functions, sharp constants, symmetry or symmetry breaking phenomenon: when $b=a+2$, then \eqref{ckn2} reduces into Rellich type inequality, Caldiroli and Musina \cite{CM12} obtained the explicit form of sharp constant when $b=a+2$ and it is not achieved. While the authors \cite{CM11} proved when $a<b<a+2$ or $a=b$ with energy assumption, the best constant is always achieved in cones, see also
    \cite[Theorem A.2]{MS14}. Furthermore, Szulkin and Waliullah \cite{SW12} proved when $a=b>0$ the best constant is achieved.
    Caldiroli and Cora \cite{CC16} obtained a partial symmetry breaking result when the parameter of pure Rellich term is sufficiently large, we refer to \cite{CM11} in cones, and also \cite{Ya21} with Hardy and Rellich terms. However, as mentioned previous, there are no optimal results about symmetry or symmetry breaking phenomenon.

    Let us recall the classical second-order Sobolev inequality
    \begin{equation*}%\label{cssi}
    \int_{\mathbb{R}^N}|\Delta u|^2 \mathrm{d}x \geq \mathcal{S}_0\left(\int_{\mathbb{R}^N}|u|^{2^{**}} \mathrm{d}x\right)^{\frac{2}{2^{**}}}, \quad \mbox{for all}\  u\in C^\infty_0(\mathbb{R}^N),
    \end{equation*}
    where $N\geq 5$, $2^{**}:=\frac{2N}{N-4}$ and $\mathcal{S}_0$ is the best constant given as
    \begin{equation}\label{cssi}
    \mathcal{S}_0=\pi^2N(N-4)(N^2-4)
    \left(\frac{\Gamma(N/2)}{\Gamma(N)}\right)^{4/N},
    \end{equation}
    see \cite{Va93}. Here $\Gamma$ denotes the usual Gamma function. It is well known that $\mathcal{S}_0$ is achieved only by $C\lambda^{\frac{N-4}{2}}
    (1+\lambda^2|x-x_0|^2)^{-\frac{N-4}{2}}$, see  \cite{Li85-1,EFJ90}. Furthermore, for the weighted Sobolev inequality (or the so-called Rellich-Sobolev inequality)
    \begin{equation*}%\label{ssi}
    \int_{\mathbb{R}^N}|\Delta u|^2 \mathrm{d}x \geq \mathcal{S}_\alpha\left(\int_{\mathbb{R}^N}
    \frac{|u|^{2^{**}_{\alpha}}}{|x|^{\alpha}} \mathrm{d}x\right)^{\frac{2}{2^{**}_{\alpha}}}, \quad \mbox{for all}\ u\in C^\infty_0(\mathbb{R}^N),
    \end{equation*}
    where $N\geq 5$, $0<\alpha<4$ and $2^{**}_{\alpha}=\frac{2(N-\alpha)}{N-4}$, it is also known that the best constant $\mathcal{S}_\alpha$ is achieved by radial functions (see the classical result of Lions \cite{Li85-2}), however the explicit form of minimizer is not known yet only for its asymptotic behavior, we refer to \cite{JL14} for details, and also \cite{DX17} with pure Rellich potentials.
    Recently, in \cite{DT23-jfa} we have considered the following second-order (CKN) type inequality:
    \begin{equation*}
    \int_{\mathbb{R}^N}|x|^{\alpha}|\Delta u|^2 \mathrm{d}x \geq S_1^{rad}(N,\alpha)\left(\int_{\mathbb{R}^N}
    |x|^{-\alpha}|u|^{p^*_{\alpha}} \mathrm{d}x\right)^{\frac{2}{p^*_{\alpha}}}, \quad \mbox{for all}\ u(x)=u(|x|)\in C^\infty_0(\mathbb{R}^N),
    \end{equation*}
    where $N\geq 3$, $4-N<\alpha<2$, $p^*_{\alpha}=\frac{2(N-\alpha)}{N-4+\alpha}$.
    As in \cite{FS03} which deals with the first-order (CKN) inequality, we restrict it in radial space, then consider the best constant $S_1^{rad}(N,\alpha)$ and its minimizer $V$ by using the change of variable $v(s)=u(s^{\frac{2}{2-\alpha}})$ which transfers it into the standard second-order Sobolev inequality
    \begin{equation*}
    \int^\infty_0\left[v''(s)+\frac{M-1}{s}v'(s)\right]^2 s^{M-1}\mathrm{d}s
    \geq \mathcal{B}(M)\left(\int^\infty_0|v(s)|^{\frac{2M}{M-4}}
    s^{M-1}\mathrm{d}s\right)^{\frac{M-4}{M}},
    \end{equation*}
    for some $\mathcal{B}(M)>0$, where $M=\frac{2N-2\alpha}{2-\alpha}>4$,
    and also classify the solutions of related linearized problem:
    \begin{equation*}
    \Delta(|x|^{\alpha}\Delta v)=(p^*_{\alpha}-1)|x|^{-\alpha} V^{p^*_{\alpha}-2}v \quad \mbox{in}\ \mathbb{R}^N,
    \end{equation*}
    and we have proved that if $\alpha$ is a negative even integer, there exist new solutions which ``replace'' the ones due to the translations invariance. Moreover, in \cite{DGT23-jde} we have also considered another case, by using the change of variable $v(s)=r^{2-N}u(r)$ with $r=s^{\frac{2}{2-\alpha}}$,
    \begin{equation*}
    \int_{\mathbb{R}^N}|x|^{\alpha}|\Delta u|^2 dx \geq S_2^{rad}(N,\alpha)\left(\int_{\mathbb{R}^N}|x|^{l}|u|^{q^*_{\alpha}} \mathrm{d}x\right)^{\frac{2}{q^*_{\alpha}}}, \quad \mbox{for all}\ u(x)=u(|x|)\in C^\infty_0(\mathbb{R}^N),
    \end{equation*}
    where $N\geq 3$, $2<\alpha<N$, $l=\frac{4(\alpha-2)(N-2)}{N-\alpha}-\alpha$ and $q^*_{\alpha}=\frac{2(N+l)}{N-4+\alpha}$. Furthermore, in \cite{DT24} we considered a new type weighted fourth-order equation, and established a new second-order (CKN) type inequality without restricting in radial space, i.e., for $N\geq 5$ and  $0<\alpha<2$,
    \begin{equation*}
    \int_{\mathbb{R}^N} |\mathrm{div}(|x|^\alpha\nabla u)|^2 \mathrm{d}x \geq C \left(\int_{\mathbb{R}^N}|u|^{\frac{2N}{N-4+2\alpha}} \mathrm{d}x\right)^{\frac{N-4+2\alpha}{N}},\quad
    \mbox{for all}\ u\in C^\infty_0(\mathbb{R}^N).
    \end{equation*}

    Therefore, it is natural to consider whether we can establish some second-order inequalities and obtain the explicit symmetry or symmetry breaking region.

\subsection{Problem setup and main results}\label{subsectmr}

    In present paper, we do not directly deal with the high order (CKN) inequality \eqref{ckn2y} established by Lin \cite{Li86}, but we will establish a new second-order (CKN) type inequality as the following.

    \begin{theorem}\label{thmgeqd}
    There is a constant $\mathcal{S}>0$ such that
    \begin{equation}\label{ckn2n}
    \int_{\mathbb{R}^N}|x|^{-\beta}|\mathrm{div} (|x|^{\alpha}\nabla u)|^2 \mathrm{d}x
    \geq \mathcal{S}\left(\int_{\mathbb{R}^N}
    |x|^{\beta}|u|^{p^{**}_{\alpha,\beta}} \mathrm{d}x\right)^{\frac{2}{p^{**}_{\alpha,\beta}}},\quad
    \mbox{for all}\ u\in C^\infty_0(\mathbb{R}^N),
    \end{equation}
    where $p^*_{\alpha,\beta}=\frac{2(N+\beta)}{N-4+2\alpha-\beta}$, and
    \begin{align}\label{cknc}
    N\geq 5,\quad \alpha>2-N,\quad \alpha-2< \beta\leq \frac{N}{N-2}\alpha.
    \end{align}
    \end{theorem}

    Note that the conditions in \eqref{cknc} are quite similar to those in Theorem A. In fact, we will show in Lemma \ref{lemgeqd} that
    \begin{align*}
    \int_{\mathbb{R}^N}|x|^{2\alpha-\beta}|\Delta u|^2 \mathrm{d}x
    \leq C\int_{\mathbb{R}^N}|x|^{-\beta}|\mathrm{div} (|x|^{\alpha}\nabla u)|^2 \mathrm{d}x,\quad \mbox{for all}\ u\in C^\infty_0(\mathbb{R}^N).
    \end{align*}
    Taking $a=-\frac{2\alpha-\beta}{2}$ and $b=-\frac{\beta}{p^{**}_{\alpha,\beta}}$, from \eqref{cknc} we have $a<\frac{N-4}{2}$, $a\leq b< a+2$ and $p^{**}_{\alpha,\beta}=\frac{2N}{N-2(2+a-b)}$, moreover, if $\alpha-2=\beta$ then $p^*_{\alpha,\beta}=2$ and $b=a+2$ which reduces \eqref{ckn2n} into Rellich type inequality.

     In order to study the sharp constant in \eqref{ckn2n}, we introduce the space $\mathcal{D}^{2,2}_{\alpha,\beta}(\mathbb{R}^N)$ as the completion of $C^\infty_0(\mathbb{R}^N)$ with respect to the norm
    \begin{equation*}
    \|u\|_{\mathcal{D}^{2,2}_{\alpha,\beta}(\mathbb{R}^N)}
    =\left(\int_{\mathbb{R}^N}|x|^{-\beta}|\mathrm{div} (|x|^{\alpha}\nabla u)|^2 \mathrm{d}x\right)^{1/2},
    \end{equation*}
    thus we rewrite \eqref{ckn2n} as
    \begin{equation}\label{ckns}
    \mathcal{S}=\inf_{u\in \mathcal{D}^{2,2}_{\alpha,\beta}(\mathbb{R}^N)\setminus\{0\}}
    \frac{\int_{\mathbb{R}^N}|x|^{-\beta}|\mathrm{div} (|x|^{\alpha}\nabla u)|^2 \mathrm{d}x}
    {\left(\int_{\mathbb{R}^N}|x|^{\beta}|u|^{p^*_{\alpha,\beta}} \mathrm{d}x\right)^{\frac{2}{p^*_{\alpha,\beta}}}}>0.
    \end{equation}
    We are interested in whether the extremals (if exist) of the best constant $\mathcal{S}$ are symmetry.

    Following the work of Felli and Schneider \cite{FS03}, firstly let us consider the radial case.
    The Euler-Lagrange equation of (CKN) inequality \eqref{ckn2n}, up to scaling, is given by
    \begin{equation}\label{Pwh}
    \mathrm{div}(|x|^{\alpha}\nabla(|x|^{-\beta}
    \mathrm{div}(|x|^\alpha\nabla u)))=|x|^\beta|u|^{p^*_{\alpha,\beta}-2}u \quad \mbox{in}\ \mathbb{R}^N,\quad u\in \mathcal{D}^{2,2}_{\alpha,\beta}(\mathbb{R}^N),
    \end{equation}
    which is equivalent to a special weighted Lane-Emden system
    \begin{eqnarray*}
    \left\{ \arraycolsep=1.5pt
       \begin{array}{ll}
        -\mathrm{div}(|x|^\alpha\nabla u)=|x|^{\beta}v\quad \mbox{in}\  \mathbb{R}^N,\\[2mm]
        -\mathrm{div}(|x|^{\alpha}\nabla v)=|x|^\beta|u|^{p^{*}_{\alpha,\beta}-2}u\quad \mbox{in}\  \mathbb{R}^N.
        \end{array}
    \right.
    \end{eqnarray*}

    \begin{theorem}\label{thmpwh}
    Assume that \eqref{cknc} holds. Then problem \eqref{Pwh} has a unique (up to scalings and change of sign)
    nontrivial radial solution of the form $\pm U_{\lambda}$ for all $\lambda>0$, where $U_{\lambda}(x)=\lambda^{\frac{N-4+2\alpha-\beta}{2}}U(\lambda x)$ with
    \begin{equation}\label{defula}
    U(x)=\frac{C_{N,\alpha,\beta}}
    {(1+|x|^{2+\beta-\alpha})
    ^{\frac{N-4+2\alpha-\beta}{2+\beta-\alpha}}}.
    \end{equation}
    Here
    $C_{N,\alpha,\beta}=\left[(N-4+2\alpha-\beta)(N-2+\alpha)
    (N+\beta)(N+2-\alpha+2\beta)\right]
    ^{\frac{N-4+2\alpha-\beta}{4(2+\beta-\alpha)}}$.
    \end{theorem}

    Therefore, as a direct consequence of Theorem \ref{thmpwh}, we obtain
    \begin{corollary}\label{thmPbcb}
    Assume that \eqref{cknc} holds. Define the best constant in the radial class as
    \begin{equation}\label{Ppbcm}
    \mathcal{S}_r=\inf_{u\in \mathcal{D}^{2,2}_{\alpha,\beta}(\mathbb{R}^N)\setminus\{0\},
    \
    u(x)=u(|x|)}
    \frac{\int_{\mathbb{R}^N}|x|^{-\beta}|\mathrm{div} (|x|^{\alpha}\nabla u)|^2 \mathrm{d}x}
    {\left(\int_{\mathbb{R}^N}|x|^{\beta}|u|^{p^*_{\alpha,\beta}} \mathrm{d}x\right)^{\frac{2}{p^*_{\alpha,\beta}}}}.
    \end{equation}
    Then the explicit form of $\mathcal{S}_r$ is
    \begin{equation*}
    \mathcal{S}_r
    =\left(\frac{2}{2+\beta-\alpha}\right)
    ^{\frac{2(2+\beta-\alpha)}{N+\beta}-4}
    \left(\frac{2\pi^{\frac{N}{2}}}{\Gamma(\frac{N}{2})}\right)
    ^{\frac{2(2+\beta-\alpha)}{N+\beta}}
    \mathcal{B}\left(\frac{2(N+\beta)}{2+\beta-\alpha}\right),
    \end{equation*}
    where $\mathcal{B}(M)=(M-4)(M-2)M(M+2)
    \left[\Gamma^2(\frac{M}{2})/(2\Gamma(M))\right]^{\frac{4}{M}}$ and $\Gamma$ is the Gamma function.
    Moreover the extremal functions which achieve $\mathcal{S}_r$ in (\ref{Ppbcm}) are given by
    \begin{equation*}%\label{pbcm}
    V_{\lambda}(x)
    =\frac{A\lambda^{\frac{N-4+2\alpha-\beta}{2}}}
    {(1+\lambda^{2+\beta-\alpha}|x|^{2+\beta-\alpha})
    ^{\frac{N-4+2\alpha-\beta}{2+\beta-\alpha}}}
    \end{equation*}
    for any $A\in\mathbb{R}\setminus\{0\}$ and $\lambda>0$.
    \end{corollary}

    Then we concern the linearized problem related to \eqref{Pwh} at the function $U$. This leads to study the problem
    \begin{equation}\label{Pwhl}
    \mathrm{div}(|x|^{\alpha}\nabla(|x|^{-\beta}
    \mathrm{div}(|x|^\alpha\nabla v)))=(p^*_{\alpha,\beta}-1)|x|^\beta U^{p^*_{\alpha,\beta}-2}v \quad \mbox{in}\ \mathbb{R}^N,\quad v\in \mathcal{D}^{2,2}_{\alpha,\beta}(\mathbb{R}^N).
    \end{equation}
    It is easy to verify that $\frac{N+2\alpha-\beta-4}{2}U+x\cdot \nabla U$ (which equals $\frac{\partial U_{\lambda}}{\partial \lambda}|_{\lambda=1}$) solves the linear equation \eqref{Pwhl}. We say $U$ is non-degenerate if all the solutions of \eqref{Pwh} result from the invariance (up to scalings) of \eqref{Pwhl}. The non-degeneracy of solutions has its own interests, such as it is a key ingredient in analyzing the blow-up phenomena of solutions to various elliptic equations on bounded or unbounded domain in $\mathbb{R}^N$ whose asymptotic behavior is encoded by $U$.
    Therefore, it is quite natural to ask the following question:
    \begin{center}
    {\em is solution $U$ non-degenerate?}
    \end{center}
    Let us define a function
    \begin{align}\label{defbfs}
    \beta_{\mathrm{FS}}(\alpha):=
        -N+\sqrt{N^2+\alpha^2+2(N-2)\alpha}.
    \end{align}
    We give an affirmative answer to above question when $\alpha-2<\beta<\beta_{\mathrm{FS}}(\alpha)$ if $\alpha>0$ and $\alpha-2< \beta\leq \frac{N}{N-2}\alpha$ if $2-N<\alpha<0$, however when $\beta=\beta_{\mathrm{FS}}(\alpha)$ there exist new solutions to the linearized problem that ``replace'' the ones due to the translations invariance.

    \begin{theorem}\label{thmpwhl}
    Assume that \eqref{cknc} holds, and $\alpha-2<\beta\leq \beta_{\mathrm{FS}}(\alpha)$ if $\alpha\geq 0$ and $\alpha-2< \beta\leq \frac{N}{N-2}\alpha$ if $2-N<\alpha<0$. If $\beta=\beta_{\mathrm{FS}}(\alpha)$, then the space of solutions of (\ref{Pwhl}) has dimension $(N+1)$ and is spanned by
    \begin{equation}\label{defaezki}
    Z_{0}(x)=\frac{1-|x|^{2+\beta-\alpha}}{(1+|x|^{2+\beta-\alpha})
    ^\frac{N-2+\alpha}{2+\beta-\alpha}},\quad Z_{i}(x)=\frac{|x|^{\frac{2+\beta-\alpha}{2}}}{(1+|x|^{2+\beta-\alpha})
    ^\frac{N-2+\alpha}{2+\beta-\alpha}}\cdot \frac{x_i}{|x|},\quad i=1,\ldots,N.
    \end{equation}
    Otherwise, the space of solutions of \eqref{Pwhl} has dimension $1$ and is spanned by $Z_0\thicksim \frac{\partial U_{\lambda}}{\partial \lambda}|_{\lambda=1}$, and in this case
    we say the solution $U$ of equation \eqref{Pwh} is non-degenerate.
    \end{theorem}

    \begin{remark}\label{rme}
    It is worth to notice that $Z_i\not\sim \frac{\partial U}{\partial x_i}$ for every $i\in \{1,\ldots,N\}$, except for $\alpha=\beta=0$. Furthermore, the curve $\beta_{\mathrm{FS}}(\alpha)$ is nothing else, it's just Felli-Schneider curve $b_{\mathrm{FS}}(a)$ as in \eqref{deffsc} by taking $\alpha=-2a$ and $\beta=-b\tau$.
    \end{remark}

    Now, we are ready to give the main result of this paper.

    \begin{theorem}\label{thmmr}
    Let $N\geq 5$, $\alpha>0$ and $\beta_{\mathrm{FS}}(\alpha)<\beta< \frac{N}{N-2}\alpha$. Then $\mathcal{S}<\mathcal{S}_r$, that is, the extremal function of the bast constant $\mathcal{S}$ defined in \eqref{ckns}, if it exists, is nonradial.
    \end{theorem}

    For the existence or non-existence of extremal functions for sharp constant $\mathcal{S}$, we can follow the work of Catrina and Wang \cite{CW01}, and also our recent work \cite{DT24-2} with some minor changes, that is, we can prove $\mathcal{S}$ is not attained if $\beta=\frac{N}{N-2}\alpha$ and $\alpha>0$, otherwise $\mathcal{S}$ is attained. Comparing our result of Theorem \ref{thmmr} to the ones in \cite{DEL16,FS03}, we also call $\beta_{\mathrm{FS}}(\alpha)$ which is given in \eqref{defbfs} as {\em Felli-Schneider curve}, then we give the following conjecture.

\vskip0.25cm

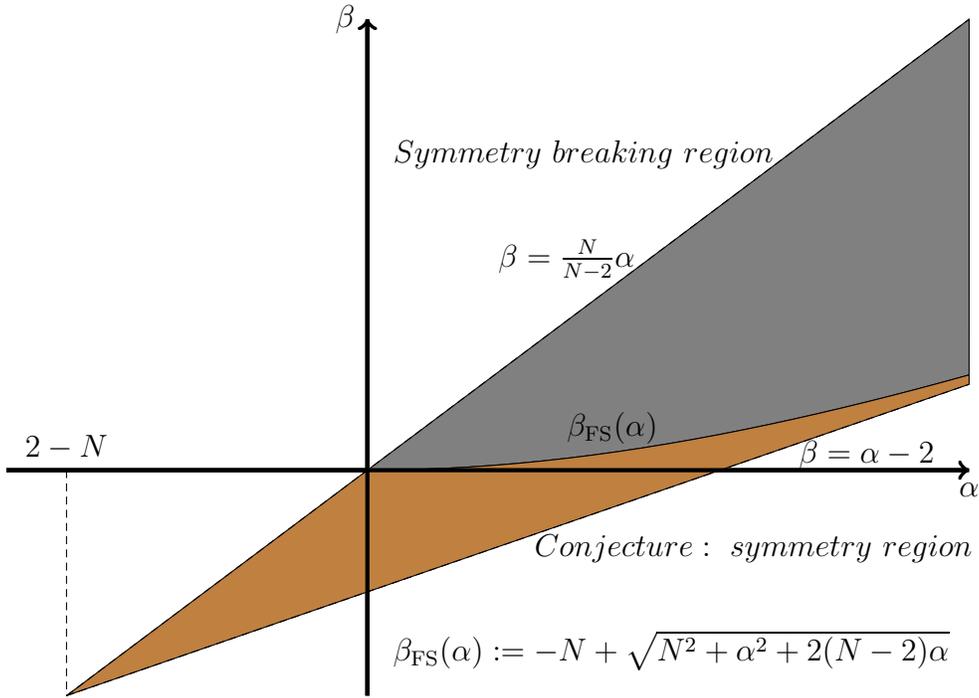
\begin{figure}[ht]
    \begin{tikzpicture}[scale=4]
%%==ÉèÖÃ×ø±êÖ᷶Χ==
		\draw[->,ultra thick](-1.2,0)--(2,0)node[below]{$\alpha$};
		\draw[->,ultra thick](0,-0.75)--(0,1.5)node[left]{$\beta$};
%%==»­Í¼==
\draw[fill=gray,domain=0:2]plot(\x,{-1+((3/7*\x)^2+1^2)^0.5})
--(2,1.5)--(0,0);
\draw[fill=brown,domain=0:2]plot(\x,{-1+((3/7*\x)^2+1^2)^0.5})
--(2,2/7)--(34/29,0);
\draw[fill=brown,domain=-1:1](-1,-0.75)--(0,0)--(34/29,0)--(-1,-0.75);
%%==Ï߶Î==
        \draw[densely dashed](-1,-3/4)--(2,1.5);
        \draw[densely dashed](-1,-3/4)--(2,2/7);
        \draw[densely dashed](-1,-0)node[above]{$2-N$}--(-1,-0.75);
        \draw[-,ultra thick](-1.2,0)--(2,0);
        \draw[-,ultra thick](0,-0.75)--(0,1.5);
%%==ÉèÖÃͼÀý==
		\node[left] at(1.92,0.05){$\beta=\alpha-2$};
        \node[right] at (0.4,0.7){$\beta=
        \frac{N}{N-2}\alpha$};
        \node[left] at (1,0.14){$\beta_{\mathrm{FS}}(\alpha)$};
        \node[right] at (0.05,-0.6){$\beta_{\mathrm{FS}}(\alpha):=
        -N+\sqrt{N^2+\alpha^2+2(N-2)\alpha}$};
        \node[right] at (0.05,1.05){$Symmetry\ breaking\ region$};
        \node[right] at (0.52,-0.26){$Conjecture:\ symmetry\ region$};
\end{tikzpicture}
\caption{\small The new second-order case, which is equivalent to Figure \ref{F1} by making the changes $\alpha=-2a$ and $\beta=-b\tau$.}
\label{F2}
\end{figure}

    {\em {\bf Conjecture:} the extremal functions of $\mathcal{S}$ are symmetry radial either $2-N<\alpha\leq 0$ and $\alpha-2<\beta\leq \frac{N}{N-2}\alpha$ (except for $\alpha=\beta=0$), or $\alpha>0$ and $\alpha-2< \beta\leq \beta_{\mathrm{FS}}(\alpha)$.} See Figure \ref{F2}.

\vskip0.25cm

    Although so far we can not provide a complete proof of this conjecture, instead, we will give a partial result about symmetry of extremal functions when $\beta= \frac{N}{N-2}\alpha$ and $\alpha<0$. Note that $\beta= \frac{N}{N-2}\alpha$ implies $p^{*}_{\alpha,\beta}=2^{**}$.

    \begin{theorem}\label{thm2ps}
    Assume that $N\geq 5$ and $2-N<\alpha<0$. For all $u\in \mathcal{D}^{2,2}_{\alpha,\frac{N}{N-2}\alpha}(\mathbb{R}^N)$,
    \begin{equation}\label{ckn2ps}
    \int_{\mathbb{R}^N}|x|^{-\frac{N}{N-2}\alpha}|\mathrm{div} (|x|^{\alpha}\nabla u)|^2 \mathrm{d}x
    \geq \left(1+\frac{\alpha}{N-2}\right)^{4-\frac{4}{N}}\mathcal{S}_0
    \left(\int_{\mathbb{R}^N}
    |x|^{\frac{N}{N-2}\alpha}|u|^{2^{**}} \mathrm{d}x\right)^{\frac{2}{2^{**}}},
    \end{equation}
    where $\mathcal{S}_0$ is the classical sharp constant of second-order Sobolev inequality given as in \eqref{cssi}. Furthermore, the constant $\left(1+\frac{\alpha}{N-2}\right)^{4-\frac{4}{N}}\mathcal{S}_0$ is sharp and equality holds if and only if $u(x)=AU_\lambda$
    for all $A\in\mathbb{R}$ and $\lambda>0$, where $U_\lambda(x)=\lambda^{\frac{N-4}{2}(1+\frac{\alpha}{N-2})}U(\lambda x)$ is as in Theorem \ref{thmpwh} replacing $\beta$ by $\frac{N}{N-2}\alpha$.
    \end{theorem}

    \begin{remark}\label{remkspnckn}\rm
    The proof of Theorem \ref{thm2ps} mainly follows from the work of Dan, Ma and Yang \cite[Theorem 1.6]{DMY20} which states that for $N\geq 5$ and $0<\mu<N-4$,
    \begin{align}\label{RSi}
    & \int_{\mathbb{R}^N}|\Delta v|^2 \mathrm{d}x
    -C_{\mu,1}\int_{\mathbb{R}^N}\frac{|\nabla v|^2}{|x|^2} \mathrm{d}x
    +C_{\mu,2}\int_{\mathbb{R}^N}\frac{|v|^2}{|x|^4} \mathrm{d}x
    \nonumber\\
    & \geq \left(1-\frac{\mu}{N-4}\right)^{4-\frac{4}{N}}\mathcal{S}_0
    \left(\int_{\mathbb{R}^N}|v|^{\frac{2N}{N-4}} \mathrm{d}x\right)^\frac{N-4}{N},\quad \mbox{for all}\ v\in \mathcal{D}^{2,2}_{0,0}(\mathbb{R}^N),
    \end{align}
    where
    \begin{align*}%\label{RSisc}
    C_{\mu,1}:& =\frac{N^2-4N+8}{2(N-4)^2}\mu[2(N-4)-\mu];
    \\
    C_{\mu,2}:& =\frac{N^2}{16(N-4)^2}\mu^2[2(N-4)-\mu]^2
    -\frac{N-2}{2}\mu[2(N-4)-\mu],
    \end{align*}
    moreover, equality in \eqref{RSi} holds if and only if
    \begin{align*}%\label{defefu}
    v(x)=A|x|^{-\frac{\mu}{2}}
    \left(\nu+|x|^{2(1-\frac{\mu}{N-4})}
    \right)^{-\frac{N-4}{2}},\quad \mbox{for}\ A\in\mathbb{R},\ \nu>0.
    \end{align*}
    The key step of proving Theorem \ref{thm2ps} is the change of variable
    \begin{align*}
    u(x)=|x|^{\eta}v(x) \quad\mbox{with}\quad
    \eta=-\frac{N-4}{2(N-2)}\alpha.
    \end{align*}
    The purpose of selecting such $\eta$ is to ensure that $2(\alpha+\eta)-\frac{N}{N-2}\alpha=\frac{N}{N-2}\alpha+\eta\cdot 2^{**}=0$.
    By some direct calculations we will prove
    \begin{align*}
    \int_{\mathbb{R}^N}|x|^{-\frac{N}{N-2}\alpha}|\mathrm{div} (|x|^{\alpha}\nabla u)|^2 \mathrm{d}x
    \geq \int_{\mathbb{R}^N}|\Delta v|^2\mathrm{d}x
    -C_{\mu,1}\int_{\mathbb{R}^N}\frac{|\nabla v|^2}{|x|^{2}} \mathrm{d}x
    +C_{\mu,2}\int_{\mathbb{R}^N}\frac{|v|^2}{|x|^{4}} \mathrm{d}x,
    \end{align*}
    and the equality holds if and only if $u$ is radially symmetry (the assumption $\alpha<0$ will play a crucial role), where $\mu=\frac{N-4}{2-N}\alpha$. Furthermore,
    \begin{align*}
    \int_{\mathbb{R}^N}
    |x|^{\frac{N}{N-2}\alpha}|u|^{2^{**}} \mathrm{d}x
    =\int_{\mathbb{R}^N}
    |v|^{2^{**}} \mathrm{d}x.
    \end{align*}
    Then our conclusion directly follows from \cite[Theorem 1.6]{DMY20}.
    \end{remark}

\subsection{Structure of this paper}\label{subsectsop}

    The paper is organized as follows:
    In Section \ref{sectnckn2} we establish the new second-order (CKN) type inequality which does not need such condition \eqref{defac}.
    In Section \ref{sectpmr} we give the uniqueness of nontrivial radial solutions of Euler-Lagrange equation \eqref{Pwh} and deduce the extremal functions of $\mathcal{S}_r$. Section \ref{sectsbp} is devoted to characterizing all solutions to the linearized problem \eqref{Pwhl}, then by using this result we will give the proof of symmetry breaking conclusion of Theorem \ref{thmmr}. Section \ref{sectps} is devoted to giving a symmetry result of extremal functions about our conjecture and proving Theorem \ref{thm2ps}.

\section{{\bfseries New second-order (CKN) type inequality}}\label{sectnckn2}

In this section, we are going to establish the new second-order (CKN) type inequality \eqref{ckn2n}. Firstly, we need the following lemma.

\begin{lemma}\label{lemgeqd}
Assume that \eqref{cknc} holds. Then there is a constant $C>0$ such that
    \begin{align}\label{neq}
    \int_{\mathbb{R}^N}|x|^{2\alpha-\beta}|\Delta u|^2 \mathrm{d}x
    \leq C\int_{\mathbb{R}^N}|x|^{-\beta}|\mathrm{div} (|x|^{\alpha}\nabla u)|^2 \mathrm{d}x,\quad \mbox{for all}\ u\in C^\infty_0(\mathbb{R}^N).
    \end{align}
\end{lemma}

\begin{proof}
We follow the arguments as those in \cite[Section 2]{GG22}. For $u\in C^\infty_0(\mathbb{R}^N)$, let
\begin{align}\label{neqer}
w:=-|x|^{-\frac{\beta}{2}}\mathrm{div} (|x|^{\alpha}\nabla u),
\end{align}
then we see that
\begin{align}\label{neqer1}
\int_{\mathbb{R}^N}|w|^2\mathrm{d}x
=\int_{\mathbb{R}^N}|x|^{-\beta}|\mathrm{div} (|x|^{\alpha}\nabla u)|^2\mathrm{d}x.
\end{align}
It follows from \eqref{neqer} that
\begin{align*}
|x|^{\alpha-\frac{\beta}{2}}\Delta u=-w-\alpha|x|^{\alpha-\frac{\beta}{2}-2}(x\cdot \nabla u),
\end{align*}
thus
\begin{align}\label{neqer2}
\int_{\mathbb{R}^N}|x|^{2\alpha-\beta}|\Delta u|^2\mathrm{d}x
& = \int_{\mathbb{R}^N}|w|^2\mathrm{d}x
+2\alpha\int_{\mathbb{R}^N}|x|^{\alpha-\frac{\beta}{2}-2}(x\cdot \nabla u)w\mathrm{d}x
\nonumber\\
& \quad +\alpha^2\int_{\mathbb{R}^N} |x|^{2\alpha-\beta-4}(x\cdot \nabla u)^2\mathrm{d}x.
\end{align}
Taking $\phi=|x|^{\alpha-\beta-2} u$ as a test function into \eqref{neqer} which is equivalent to $|x|^{\frac{\beta}{2}}w=-\mathrm{div} (|x|^{\alpha}\nabla u)$,
\[
\int_{\mathbb{R}^N}|x|^{\frac{\beta}{2}}w \phi\mathrm{d}x
=\int_{\mathbb{R}^N}|x|^{\alpha}\nabla u\cdot \nabla \phi\mathrm{d}x,
\]
therefore,
\begin{align}\label{neqere}
\int_{\mathbb{R}^N} |x|^{\alpha-\frac{\beta}{2}-2}wu\mathrm{d}x
& =\int_{\mathbb{R}^N} |x|^{2\alpha-\beta-2}|\nabla u|^2\mathrm{d}x+(\alpha-\beta-2)\int_{\mathbb{R}^N}|x|^{2\alpha-\beta-4}
\left(x\cdot \nabla \left(\frac{|u|^2}{2}\right)\right)\mathrm{d}x
\nonumber\\
& =\frac{(2+\beta-\alpha)(N+2\alpha-\beta-4)}{2}
\int_{\mathbb{R}^N}|x|^{2\alpha-\beta-4}|u|^2\mathrm{d}x
\nonumber\\
& \quad +\int_{\mathbb{R}^N} |x|^{2\alpha-\beta-2}|\nabla u|^2\mathrm{d}x.
\end{align}
Note that the assumptions imply $\frac{(2+\beta-\alpha)(N+2\alpha-\beta-4)}{2}\geq 0$,
then it follows from \eqref{neqere} and the H\"{o}lder inequality that
\begin{align*}
\int_{\mathbb{R}^N} |x|^{2\alpha-\beta-2}|\nabla u|^2\mathrm{d}x
\leq \int_{\mathbb{R}^N} |x|^{\alpha-\frac{\beta}{2}-2}wu\mathrm{d}x
\leq \left(\int_{\mathbb{R}^N} |w|^2\mathrm{d}x\right)^{\frac{1}{2}}\left(\int_{\mathbb{R}^N} |x|^{2\alpha-\beta-4}|u|^2\mathrm{d}x\right)^{\frac{1}{2}}.
\end{align*}
Since the assumptions also imply $-\frac{2\alpha-\beta-2}{2}<\frac{N-2}{2}$, by using the weighted Hardy inequality (see \cite{CW01}) we obtain
\[
\int_{\mathbb{R}^N}|x|^{2\alpha-\beta-4}|u|^2\mathrm{d}x
\leq E \int_{\mathbb{R}^N}|x|^{2\alpha-\beta-2}|\nabla u|^2\mathrm{d}x,
\]
for some $E>0$ independent of $u$ (in fact, from \cite{CW01} we know the sharp constant is $E=(\frac{2}{N-4+2\alpha-\beta})^2$), then we see that
\begin{align*}
\int_{\mathbb{R}^N} |x|^{2\alpha-\beta-2}|\nabla u|^2\mathrm{d}x
\leq E\int_{\mathbb{R}^N} |w|^2\mathrm{d}x.
\end{align*}
By the Young inequality and $(x\cdot\nabla u)^2\leq |x|^2|\nabla u|^2$, it follows from \eqref{neqer1} and \eqref{neqer2} that
\begin{align*}
\int_{\mathbb{R}^N}|x|^{2\alpha-\beta}|\Delta u|^2\mathrm{d}x
& \leq \int_{\mathbb{R}^N}|w|^2\mathrm{d}x
+|\alpha|\left(\int_{\mathbb{R}^N}|x|^{2\alpha-\beta-2}|\nabla u|^2\mathrm{d}x
+\int_{\mathbb{R}^N}|w|^2\mathrm{d}x
\right)
\nonumber\\
& \quad +\alpha^2\int_{\mathbb{R}^N} |x|^{2\alpha-\beta-2}|\nabla u|^2\mathrm{d}x
\\
& \leq \left[1+|\alpha|+E(|\alpha|+\alpha^2)\right]\int_{\mathbb{R}^N}
|x|^{-\beta}|\mathrm{div} (|x|^{\alpha}\nabla u)|^2\mathrm{d}x.
\end{align*}
Therefore, we obtain inequality \eqref{neq} with $C=1+|\alpha|+E(|\alpha|+\alpha^2)>0$.
\end{proof}

From Lemma \ref{lemgeqd} and second-order Caffarelli-Kohn-Nirenberg inequality \eqref{ckn2}, if condition \eqref{defac} with $2\alpha-\beta=-2a$ holds, then the new second-order Caffarelli-Kohn-Nirenberg type inequality \eqref{ckn2n} follows.
Next, we will show \eqref{ckn2n} holds which does not need the additional condition \eqref{defac} with $2\alpha-\beta=-2a$.

\subsection{Proof of Theorem \ref{thmgeqd}}\label{subsectpeckn2}
When $2-N<\alpha\leq 0$, it is obvious that $-\frac{N}{2}<-\frac{2\alpha-\beta}{2}<\frac{N-4}{2}$, thus the condition \eqref{defac} always holds with $a=-\frac{2\alpha-\beta}{2}$.
Combining with \eqref{neq} and the second-order CKN inequality \eqref{ckn2} so that \eqref{ckn2n} holds.

Now, we only need to consider $\alpha>0$. By a direct calculation we have $\mathrm{div} (|x|^{\alpha}\nabla u)=|x|^\alpha\Delta u+\alpha |x|^{\alpha-2}(x\cdot\nabla u)$, then
    \begin{align}\label{ckn2nfe}
    \int_{\mathbb{R}^N}|x|^{-\beta}|\mathrm{div} (|x|^{\alpha}\nabla u)|^2\mathrm{d}x
    & = \int_{\mathbb{R}^N}|x|^{2\alpha-\beta}|\Delta u|^2\mathrm{d}x
    + 2\alpha \int_{\mathbb{R}^N}|x|^{2\alpha-2-\beta}\Delta u(x\cdot\nabla u)\mathrm{d}x
    \nonumber \\
    & \quad + \alpha^2 \int_{\mathbb{R}^N}|x|^{2\alpha-4-\beta}(x\cdot\nabla u)^2\mathrm{d}x\nonumber \\
    & = \int_{\mathbb{R}^N}|x|^{2\alpha-\beta}|\Delta u|^2\mathrm{d}x
    \nonumber \\
    & \quad+ \alpha (N-4+2\alpha-\beta)\int_{\mathbb{R}^N}
    |x|^{2\alpha-2-\beta}|\nabla u|^2\mathrm{d}x
    \nonumber \\
    & \quad+ \alpha(2\beta-3\alpha+4) \int_{\mathbb{R}^N}|x|^{2\alpha-4-\beta}(x\cdot\nabla u)^2\mathrm{d}x,
    \end{align}
    thanks to
    \begin{align*}
    \int_{\mathbb{R}^N}|x|^{2\alpha-2-\beta}\Delta u(x\cdot\nabla u)\mathrm{d}x
    & = \frac{N-4+2\alpha-\beta}{2}\int_{\mathbb{R}^N}
    |x|^{2\alpha-2-\beta}|\nabla u|^2\mathrm{d}x
    \\
    & \quad + (\beta-2\alpha+2) \int_{\mathbb{R}^N}|x|^{2\alpha-4-\beta}(x\cdot\nabla u)^2\mathrm{d}x.
    \end{align*}
    Since $N-4+2\alpha-\beta\geq N-4+2\alpha-\frac{N\alpha}{N-2}=\frac{(N-4)(N-2+\alpha)}{N-2}>0$, if $2\beta-3\alpha+4\geq 0$,
     \begin{align*}
    \int_{\mathbb{R}^N}|x|^{-\beta}|\mathrm{div} (|x|^{\alpha}\nabla u)|^2\mathrm{d}x
    & \geq \int_{\mathbb{R}^N}|x|^{2\alpha-\beta}|\Delta u|^2\mathrm{d}x
    \\
    & \quad+ \alpha (N-4+2\alpha-\beta)\int_{\mathbb{R}^N}
    |x|^{2\alpha-2-\beta}|\nabla u|^2\mathrm{d}x,
    \end{align*}
    then we can deduce \eqref{ckn2n} from the first-order CKN inequality \eqref{cknwit} due to $-\frac{2\alpha-2-\beta}{2}< \frac{N-2}{2}$.

    It remains only the case $\alpha>0$ and $2\beta-3\alpha+4<0$. Thanks to $(x\cdot\nabla u)^2\leq |x|^2|\nabla u|^2$, then it follows from \eqref{ckn2nfe} that
    \begin{align*}
    \int_{\mathbb{R}^N}|x|^{-\beta}|\mathrm{div} (|x|^{\alpha}\nabla u)|^2\mathrm{d}x
    \geq & \int_{\mathbb{R}^N}|x|^{2\alpha-\beta}|\Delta u|^2\mathrm{d}x
    + \alpha (N-\alpha+\beta)\int_{\mathbb{R}^N}
    |x|^{2\alpha-2-\beta}|\nabla u|^2\mathrm{d}x.
    \end{align*}
    In this case, we have $N-\alpha+\beta\geq N-\alpha+\alpha-2=N-2>0$, then the first-order CKN inequality \eqref{cknwit} also indicates \eqref{ckn2n} holds. Now, the proof of Theorem \ref{thmgeqd} is completed.
\qed

\section{{\bfseries Extremal functions in radial case}}\label{sectpmr}

    In this section, we use a suitable transform that is changing the variable $r\mapsto r^{\frac{2}{2+\beta-\alpha}}$ related to radial Sobolev inequality, to investigate the sharp constant $\mathcal{S}_r$ and its extremal functions. Firstly, we show the uniqueness of nontrivial radial solutions of Euler-Lagrange equation \eqref{Pwh}.

\subsection{Proof of Theorem \ref{thmpwh}}\label{subsectmrs}
     Let $u\in \mathcal{\mathcal{D}}^{2,2}_{\alpha,\beta}(\mathbb{R}^N)$ be a nontrivial radial solution of equation (\ref{Pwh}). Set
    \begin{equation*}%\label{PpwhlW}
    \begin{split}
    -{\rm div}(|x|^{\alpha}\nabla u)=|x|^{\beta}v,
    \end{split}
    \end{equation*}
    and $r=|x|$, then equation (\ref{Pwh}) is equivalent to the following system:
    \begin{eqnarray}\label{Pwhlpe}
    \left\{ \arraycolsep=1.5pt
       \begin{array}{ll}
        u''(r)+\frac{N-1+\alpha}{r}u'(r)
        +\frac{v(r)}{r^{\alpha-\beta}}=0\quad \mbox{in}\ r\in(0,\infty),\\[3mm]
        v''(r)+\frac{N-1+\alpha}{r}v'(r)+ r^\beta |u|^{p^*_{\alpha,\beta}-2}u=0\quad \mbox{in}\ r\in(0,\infty).
        \end{array}
    \right.
    \end{eqnarray}
    Moreover, let us set $r=s^q$ with $q=2/(2+\beta-\alpha)$ and let
    \begin{equation*}%\label{p2txye}
    X(s)=u(r),\quad Y(s)=q^{2}v(r),
    \end{equation*}
    which transforms (\ref{Pwhlpe}) into the system
    \begin{eqnarray*}%\label{p2te}
    \left\{ \arraycolsep=1.5pt
       \begin{array}{ll}
        X''(s)+\frac{M-1}{s}X'(s)+Y(s)=0 \quad \mbox{in}\ s\in(0,\infty),\\[3mm]
        Y''(s)+\frac{M-1}{s}Y'(s)
        +q^4 |X|^{p^*_{\alpha,\beta}-2}X=0\quad \mbox{in}\ s\in(0,\infty),
        \end{array}
    \right.
    \end{eqnarray*}
    where $M=\frac{2(N+\beta)}{2+\beta-\alpha}>4$, which is equivalent to
    \begin{align}\label{PpwhlWe}
    & X^{(4)}(s)+\frac{2(M-1)}{s}X'''(s)+\frac{(M-1)(M-3)}{s^2}X''(s)
    -\frac{(M-1)(M-3)}{s^3}X'(s)
    \nonumber \\
    & \quad = q^4 |X|^{\frac{8}{M-4}}X\quad \mbox{in}\ s\in(0,\infty),\quad X\in \mathcal{C},
    \end{align}
    due to $p^*_{\alpha,\beta}=\frac{2M}{M-4}$, where $\mathcal{C}$ is defined by
    \[
    \mathcal{C}:=\left\{\omega\in C^\infty_0([0,\infty))| \int^\infty_0 \left[\omega''(s)+\frac{M-1}{s}\omega'(s)
    \right]^2
    s^{M-1}\mathrm{d}s<\infty\right\}.
    \]
    Finally, let us make the classical Emden-Fowler transformation, that is,
    \[
    X(s)=q^{-\frac{M-4}{2}}s^{-\frac{M-4}{2}}\varphi(t)
    \quad\mbox{with}\ t=-\ln s,
    \]
    which transforms \eqref{PpwhlWe} into the following fourth-order ordinary differential equation
    \begin{equation}\label{Pwht}
    \varphi^{(4)}-\frac{(M-2)^2+4}{2}\varphi''
    +\frac{M^2(M-4)^2}{16}\varphi
    =|\varphi|^{\frac{8}{M-4}}\varphi \quad \mbox{in}\  \mathbb{R},\quad \varphi\in H^2(\mathbb{R}),
    \end{equation}
    where $H^2(\mathbb{R})$ denotes the completion of $C^\infty_0(\mathbb{R})$ with respect to the norm
    \[
    \|\varphi\|_{H^2(\mathbb{R})}
    =\left[\int_{\mathbb{R}}\left(|\varphi''|^2
    +\frac{(M-2)^2+4}{2}|\varphi'|^2
    +\frac{M^2(M-4)^2}{16}|\varphi|^2\right)\mathrm{d}t\right]^{1/2}.
    \]
    Then applying
    \cite[Theorem 2.2]{BM12} directly, we could get the existence and uniqueness (up to translations, inversion $t\mapsto -t$ and change of sign) of smooth solution to equation \eqref{Pwht}, which indicates equation \eqref{PpwhlWe} has a unique (up to scalings and change of sign) nontrivial radial solution. It is easy to verify that \eqref{PpwhlWe} admits solutions of the form
    \[
    X(s)=\frac{\left[q^{-4}(M-4)(M-2)M(M+2)\right]^{\frac{M-4}{8}}\nu^{\frac{M-4}{2}}}{(1+\nu^{2}s^{2})^{\frac{M-4}{2}}}
    \]
    for all $\nu>0$.
    Thus, we deduce that the Euler-Lagrange equation \eqref{Pwh} has a unique (up to scalings and change of sign) nontrivial radial solution $\bar{u}$ of the form
    \begin{equation*}
    \bar{u}=\pm\frac{C_{N,\alpha,\beta}\lambda^{\frac{N-4+2\alpha-\beta}{2}}}
    {(1+\lambda^{2+\beta-\alpha}|x|^{2+\beta-\alpha})
    ^{\frac{N-4+2\alpha-\beta}{2+\beta-\alpha}}}
    \end{equation*}
    for all $\lambda>0$, where
    \[
    C_{N,\alpha,\beta}=\left[(N-4+2\alpha-\beta)(N-2+\alpha)
    (N+\beta)(N+2-\alpha+2\beta)\right]
    ^{\frac{N-4+2\alpha-\beta}{4(2+\beta-\alpha)}}.
    \]
    Now, the proof of Theorem \ref{thmpwh} is completed.
    \qed

\vskip0.25cm

    Then, based on the uniqueness result of Theorem \ref{thmpwh}, let us investigate the sharp constant $\mathcal{S}_r$ and its extremal functions.

\subsection{Proof of Corollary \ref{thmPbcb}}\label{subsectmrf} Let $u\in \mathcal{\mathcal{D}}^{2,2}_{\alpha,\beta}(\mathbb{R}^N)\setminus\{0\}$ be radial. As in the proof of Theorem \ref{thmpwh}, let us make the change of variable that $v(s)=u(r)$, $|x|=r=s^{q}$ with $q=2/(2+\beta-\alpha)$, then we have
    \begin{align*}
    & \int_{\mathbb{R}^N}|x|^{-\beta}|\mathrm{div} (|x|^{\alpha}\nabla u)|^2 \mathrm{d}x
    \\
    & = \omega_{N-1}\int^\infty_0 \left[ u''(r)+\frac{N+\alpha-1}{r}u'(r)\right]^2r^{N+2\alpha-\beta-1}
    \mathrm{d}r \\
    & = \omega_{N-1}q^{-3}\int^\infty_0\left[v''(s)
    +\frac{q(N-2+\alpha)+1}{s}v'(s)\right]^2
    s^{q(N-4+2\alpha-\beta)+3}\mathrm{d}s,
    \end{align*}
    and
    \begin{align*}
    \int_{\mathbb{R}^N}|x|^{\beta}|u|^{p^*_{\alpha,\beta}} \mathrm{d}x
    =\omega_{N-1}q\int^\infty_0
    |v(s)|^{p^*_{\alpha,\beta}}s^{q(N+\beta)-1}
    \mathrm{d}s.
    \end{align*}
    Here $\omega_{N-1}=2\pi^{\frac{N}{2}}/\Gamma(\frac{N}{2})$ is the surface area for unit ball of $\mathbb{R}^N$. By the chosen of $q$, it is easy to verify that
    \[
    q(N-2+\alpha)+1=q(N-4+2\alpha-\beta)+3=q(N+\beta)-1
    =\frac{2(N+\beta)}{2+\beta-\alpha}-1.
    \]
    Now, we set
    \begin{equation}\label{defm}
    M=\frac{2(N+\beta)}{2+\beta-\alpha},
    \end{equation}
    which implies $p^*_{\alpha,\beta}=\frac{2(N+\beta)}{N-4+2\alpha-\beta}
    =\frac{2M}{M-4}$, and $M>4$ due to the assumption \eqref{cknc}.
    To sum up, we deduce that
    \begin{equation}\label{defbcst}
    \mathcal{S}_r=q^{\frac{4}{M}-4}\omega^{\frac{4}{M}}_{N-1}
    \mathcal{B}(M),
    \end{equation}
    where
    \begin{equation}\label{defbcscm}
    \mathcal{B}(M)=\inf_{v\in \mathcal{\mathcal{D}}^{2,2}_{\infty}(M)\backslash\{0\}}
    \frac{\int^\infty_0\left[v''(s)+\frac{M-1}{s}v'(s)\right]^2
    s^{M-1}\mathrm{d}s}
    {\left(\int^\infty_0|v(s)|^{\frac{2M}{M-4}}s^{M-1}
    \mathrm{d}s \right)^{\frac{M-4}{M}}}.
    \end{equation}
    Here $\mathcal{\mathcal{D}}^{2,2}_{\infty}(M)$ denotes the completion of $C^\infty_0([0,\infty))$ with respect to the norm
    \[
    \|v\|^2_{\mathcal{\mathcal{D}}^{2,2}_{\infty}(M)}=\int^\infty_0\left[v''(s)+\frac{M-1}{s}v'(s)\right]^2
    s^{M-1}\mathrm{d}s.
    \]
    Following the arguments as those in our recent work \cite{DT24}, the proof can be done. In fact, de Oliveira and Silva \cite{dS24} proved that $\mathcal{B}(M)>0$ can be achieved, and the minimizers are solutions (up to some suitable multiplications) of
    \begin{equation}\label{les}
    \Delta_s^2v=|v|^{\frac{8}{M-4}}v\quad \mbox{in}\ (0,\infty),\quad v\in \mathcal{\mathcal{D}}^{2,2}_{\infty}(M),
    \end{equation}
    where $\Delta_s=\frac{\partial^2}{\partial s^2}+\frac{M-1}{s}\frac{\partial}{\partial s}$. As in the proof of Theorem \ref{thmpwh}, the uniqueness of nontrivial solutions for \eqref{Pwht} indicates that \eqref{les} admits a unique (up to scalings and change of sign) nontrivial solution of the form
    \[
    \bar{v}(s)=C\nu^{\frac{M-4}{2}}(1+\nu^2s^2)^{-\frac{M-4}{2}}, \quad \mbox{for some suitable}\ C\in\mathbb{R}\setminus\{0\}\ \mbox{and for all}\ \nu>0,
    \]
    which implies $\bar{v}$ is the unique (up to scalings and multiplications) extremal function for $\mathcal{B}(M)$. Therefore, putting $\bar{v}$ into \eqref{defbcscm} as a test function, we can directly obtain
    \begin{align*}
    &\mathcal{B}(M)=(M-4)(M-2)M(M+2)
    \left(\frac{\Gamma^2\left(\frac{M}{2}\right)}
    {2\Gamma(M)}\right)^{\frac{4}{M}},
    \end{align*}
    where $\Gamma$ denotes the classical Gamma function.
    Then turning back to \eqref{defbcst}, we have
    \begin{equation*}
    \begin{split}
    \mathcal{S}_r
    =\left(\frac{2}{2+\beta-\alpha}\right)
    ^{\frac{2(2+\beta-\alpha)}{N+\beta}-4}
    \left(\frac{2\pi^{\frac{N}{2}}}{\Gamma(\frac{N}{2})}\right)
    ^{\frac{2(2+\beta-\alpha)}{N+\beta}}
    \mathcal{B}\left(\frac{2(N+\beta)}{2+\beta-\alpha}\right),
    \end{split}
    \end{equation*}
    and it is achieved if and only if by
    \begin{equation*}
    V_{\lambda}(x)
    =\frac{A\lambda^{\frac{N-4+2\alpha-\beta}{2}}}
    {(1+\lambda^{2+\beta-\alpha}|x|^{2+\beta-\alpha})
    ^{\frac{N-4+2\alpha-\beta}{2+\beta-\alpha}}},
    \end{equation*}
    for all $A\in\mathbb{R}\backslash\{0\}$ and $\lambda>0$. The proof of Corollary \ref{thmPbcb} is now completed.
    \qed

\section{{\bfseries Symmetry breaking phenomenon}}\label{sectsbp}

In order to show the symmetry breaking phenomenon of extremal functions as Theorem \ref{thmmr}, firstly, by using the standard spherical decomposition and taking the change of variable $v(s)=u(r)$ with $r=s^{\frac{2}{2+\beta-\alpha}}$, we characterize all solutions to the linearized problem (\ref{Pwhl}).

\subsection{Proof of Theorem \ref{thmpwhl}}\label{subsectpnd}
    We follow the arguments as those in \cite{BWW03}, and also our recent work \cite{DT24}.
    Let us decompose the fourth-order equation (\ref{Pwhl}) into a system of two second-order equations. For $v\in\mathcal{D}^{2,2}_{\alpha,\beta}(\mathbb{R}^N)$, let
    \begin{align*}
    -\mathrm{div}(|x|^{\alpha}\nabla v)=|x|^\beta w,
    \end{align*}
    then problem (\ref{Pwhl}) is equivalent to the following system:
    \begin{eqnarray}\label{Pwhlp}
    \left\{ \arraycolsep=1.5pt
       \begin{array}{ll}
        -|x|^{\alpha}\Delta v-\alpha|x|^{\alpha-2}(x\cdot\nabla v)=|x|^\beta w \quad \mbox{in}\  \mathbb{R}^N,\\[3mm]
        -|x|^{\alpha}\Delta w-\alpha|x|^{\alpha-2}(x\cdot\nabla w)
        =\frac{(p^*_{\alpha,\beta}-1)C_{N,\alpha,\beta}
        ^{p^*_{\alpha,\beta}-2}}{(1+|x|^{2+\beta-\alpha})^4}
        |x|^\beta v\quad \mbox{in}\  \mathbb{R}^N.
        \end{array}
    \right.
    \end{eqnarray}

    Firstly, we decompose $v$ and $w$ as follows:
    \begin{equation}\label{defvd}
    v(r,\theta)=\sum^{\infty}_{k=0}\sum^{l_k}_{i=1}
    \phi_{k,i}(r)\Psi_{k,i}(\theta),\quad w(r,\theta)=\sum^{\infty}_{k=0}\sum^{l_k}_{i=1}
    \psi_{k,i}(r)\Psi_{k,i}(\theta),
    \end{equation}
    where $r=|x|$, $\theta=x/|x|\in \mathbb{S}^{N-1}$, and
    \begin{equation*}
    \phi_{k,i}(r)=\int_{\mathbb{S}^{N-1}}
    v(r,\theta)\Psi_{k,i}(\theta)\mathrm{d}\theta,\quad \psi_{k,i}(r)=\int_{\mathbb{S}^{N-1}}
    w(r,\theta)\Psi_{k,i}(\theta)\mathrm{d}\theta.
    \end{equation*}
    Here $\Psi_{k,i}(\theta)$ denotes the $k$-th spherical harmonic, i.e., it satisfies
    \begin{equation}\label{deflk}
    -\Delta_{\mathbb{S}^{N-1}}\Psi_{k,i}=\lambda_k \Psi_{k,i},
    \end{equation}
    where $\Delta_{\mathbb{S}^{N-1}}$ is the Laplace-Beltrami operator on $\mathbb{S}^{N-1}$ with the standard metric and  $\lambda_k$ is the $k$-th eigenvalue of $-\Delta_{\mathbb{S}^{N-1}}$. It is well known that \begin{equation}\label{deflklk}
    \lambda_k=k(N-2+k),\quad k=0,1,2,\ldots,
    \end{equation}
    whose multiplicity is
    \begin{equation}\label{deflkm}
    l_k:=\frac{(N+2k-2)(N+k-3)!}{(N-2)!k!}
    \end{equation}
    (note that $l_0:=1$) and that \[
    \mathrm{Ker}(\Delta_{\mathbb{S}^{N-1}}+\lambda_k)
    =\mathbb{Y}_k(\mathbb{R}^N)|_{\mathbb{S}^{N-1}},
    \]
    where $\mathbb{Y}_k(\mathbb{R}^N)$ is the space of all homogeneous harmonic polynomials of degree $k$ in $\mathbb{R}^N$. See \cite[Corollary C.I.3]{BGM71} for details. It is standard that $\lambda_0=0$ and the corresponding eigenfunction of (\ref{deflk}) is the constant function that is $\Psi_{0,1}=c\in\mathbb{R}\setminus\{0\}$. The second eigenvalue $\lambda_1=N-1$ and the corresponding eigenfunctions of (\ref{deflk}) are $\Psi_{1,i}=x_i/|x|$, $i=1,\ldots,N$.

    It is known that
    \begin{align}\label{Ppwhl2deflklw}
    \Delta (\varphi_{k,i}(r)\Psi_{k,i}(\theta))
    & = \Psi_{k,i}\left(\varphi''_{k,i}+\frac{N-1}{r}\varphi'_{k,i}\right)
    +\frac{\varphi_{k,i}}{r^2}\Delta_{\mathbb{S}^{N-1}}\Psi_{k,i} \nonumber\\
    & = \Psi_{k,i}\left(\varphi''_{k,i}+\frac{N-1}{r}\varphi'_{k,i}
    -\frac{\lambda_k}{r^2}\varphi_{k,i}\right).
    \end{align}
    Furthermore, it is easy to verify that
    \begin{equation*}
    \frac{\partial (\varphi_{k,i}(r)\Psi_{k,i}(\theta))}{\partial x_j}=\varphi'_{k,i}\frac{x_j}{r}\Psi_{k,i}
    +\varphi_{k,i}\frac{\partial\Psi_{k,i}}{\partial \theta_l}\frac{\partial\theta_l}{\partial x_j},\quad \mbox{for all}\quad l=1,\ldots,N-1,
    \end{equation*}
    and
    \begin{equation*}
    \sum^{N}_{j=1}\frac{\partial\theta_l}{\partial x_j}x_j=0,\quad \mbox{for all}\quad l=1,\ldots,N-1,
    \end{equation*}
    hence
    \begin{align}\label{Ppwhl2deflkln}
    x\cdot\nabla (\varphi_{k,i}(r)\Psi_{k,i}(\theta))=
    \sum^{N}_{j=1}x_j\frac{\partial (\varphi_{k,i}(r)\Psi_{k,i}(\theta))}{\partial x_j}=\varphi'_{k,i}r\Psi_{k,i}.
    \end{align}
    Therefore, by standard regularity theory, putting together (\ref{Ppwhl2deflklw}) and (\ref{Ppwhl2deflkln}) into (\ref{Pwhlp}), the function $(v,w)$ is a solution of (\ref{Pwhlp}) if and only if for all $i=1,\ldots,l_k$, the functions pair $(\phi_{k,i},\psi_{k,i})\in \mathcal{C}_k\times \mathcal{C}_k$ is a classical solution of the system
    \begin{eqnarray}\label{p2c}
    \left\{ \arraycolsep=1.5pt
       \begin{array}{ll}
        \phi''_{k,i}+\frac{N-1+\alpha}{r}\phi'_{k,i}
        -\frac{\lambda_k}{r^2}\phi_{k,i}
        +\frac{\psi_{k,i}}{r^{\alpha-\beta}}=0 \quad \mbox{in}\quad r\in(0,\infty),\\[3mm]
        \psi''_{k,i}+\frac{N-1+\alpha}{r}\psi'_{k,i}
        -\frac{\lambda_k}{r^2}\psi_{k,i}
        +\frac{(p^*_{\alpha,\beta}-1)C_{N,\alpha,\beta}
        ^{p^*_{\alpha,\beta}-2}}
        {r^{\alpha-\beta} (1+r^{2+\beta-\alpha})^4}\phi_{k,i}=0 \quad \mbox{in}\quad r\in(0,\infty),\\[3mm]
        \phi'_{k,i}(0)=\psi'_{k,i}(0)=0 \quad\mbox{if}\quad k=0,\quad \mbox{and}\quad \phi_{k,i}(0)=\psi_{k,i}(0)=0 \quad\mbox{if}\quad k\geq 1,
        \end{array}
    \right.
    \end{eqnarray}
    where $\mathcal{C}_k$ is defined by
    \[
    \mathcal{C}_k:=\left\{\omega\in C^\infty_0([0,\infty))| \int^\infty_0 \left[\omega''(r)+\frac{N+\alpha-1}{r}\omega'(r)
    -\frac{\lambda_k}{r^2}\omega(r)\right]^2
    r^{N+2\alpha-\beta-1}\mathrm{d}r<\infty\right\}.
    \]
    Then let us make the same change of variable as in the proof of Theorem \ref{thmpwh}, $|x|=r=s^q$ with $q=2/(2+\beta-\alpha)$ and let
    \begin{equation}\label{p2txy}
    X_{k,i}(s)=\phi_{k,i}(r),\quad Y_{k,i}(s)=q^2\psi_{k,i}(r),
    \end{equation}
    which transforms (\ref{p2c}) into the system
    \begin{eqnarray}\label{p2t}
    \left\{ \arraycolsep=1.5pt
       \begin{array}{ll}
        X''_{k,i}+\frac{M-1}{s}X'_{k,i}
        -\frac{q^2\lambda_k}{s^2}X_{k,i}+Y_{k,i}=0 \quad \mbox{in}\quad s\in(0,\infty),\\[3mm]
        Y''_{k,i}+\frac{M-1}{s}Y'_{k,i}
        -\frac{q^2\lambda_k}{s^2}Y_{k,i}
        +\frac{(M+4)(M-2)M(M+2)}{(1+s^2)^4}X_{k,i}=0 \quad \mbox{in}\quad s\in(0,\infty),\\[3mm]
        X'_{k,i}(0)=Y'_{k,i}(0)=0 \quad\mbox{if}\quad k=0,\quad \mbox{and}\quad X_{k,i}(0)=Y_{k,i}(0)=0 \quad\mbox{if}\quad k\geq 1,
        \end{array}
    \right.
    \end{eqnarray}
    for all $i=1,\ldots,l_k$, in $(X_{k,i},Y_{k,i})\in \widetilde{\mathcal{C}}_k\times \widetilde{\mathcal{C}}_k$, where $M=\frac{2(N+\beta)}{2+\beta-\alpha}>4$ and
    \[
    \widetilde{\mathcal{C}}_k:=\left\{\omega\in C^\infty_0([0,\infty))| \int^\infty_0 \left[\omega''(s)+\frac{M-1}{s}\omega'(s)
    -\frac{q^2\lambda_k}{s^2}\omega(s)\right]^2 s^{M-1} \mathrm{d}s<\infty\right\}.
    \]
    Here we have used the fact
    \begin{small}\begin{equation*}
    q^4(p^*_{\alpha,\beta}-1)
    C_{N,\alpha,\beta}^{p^*_{\alpha,\beta}-2}
    =\left[(M-4)(M-2)M(M+2)\right]\left[\frac{2M}{M-4}-1\right]
    =(M+4)(M-2)M(M+2).
    \end{equation*}\end{small}

    Note that \eqref{p2t} is equivalent to the fourth-order ODE
    \begin{align}\label{rwevpb}
    \left(\Delta_s-\frac{\varpi_k}{s^2}\right)^2X_{k,i}
    & = \left(q^2\lambda_k-\varpi_k\right)
    \left(\frac{2}{s^2}X_{k,i}''
    +\frac{2(M-3)}{s^3}X_{k,i}'
    -\frac{2(M-4)+q^2\lambda_k+\varpi_k}{s^4}X_{k,i}\right)
    \nonumber \\
    & \quad +(\tilde{2}^{**}-1)\Gamma_M(1+s^2)^{-4}X_{k,i},
    \end{align}
    in $s\in(0,\infty)$, $X_{k,i}\in \widetilde{\mathcal{C}}_k$, $i=1,2,\ldots, l_k$. Here $\Delta_s:=\frac{\partial^2}{\partial s^2}+\frac{M-1}{s}
    \frac{\partial}{\partial s}$, $\varpi_k:=k(M-2+k)$, $\tilde{2}^{**}:=\frac{2M}{M-4}$ and
    \begin{align}\label{defgn}
    \Gamma_M:=(M-4)(M-2)M(M+2).
    \end{align}
    Note that $q^2\lambda_k\geq\varpi_k$ for all $k\geq 1$, furthermore, $q^2\lambda_k>\varpi_k$ for $\beta<\beta_{\mathrm{FS}}(\alpha)$ and $k\geq 1$, also for $\beta=\beta_{\mathrm{FS}}(\alpha)$ and $k\geq 2$.
    It is easy to verify that when $k=0$, \eqref{rwevpb} admits only one solution $X_0(s)=\frac{1-s^2}{(1+s^2)^{\frac{M-2}{2}}}$ (up to multiplications), see our recent work \cite{DT24} for details. In fact, \cite[Lemma 2.4]{BWW03} states that if $X\in \widetilde{\mathcal{C}_0}$ is a solution of
    \[
    \left[s^{1-M}\frac{\partial}{\partial s}\left(s^{M-1}\frac{\partial}{\partial s}\right)\right]^2X=\nu(1+s^2)^{-4}X \quad \mbox{for}\ \nu>0,
    \]
    with $X(0)=0$, then $X\equiv 0$ (which states that $M$ is an integer, in fact, it also holds for all $M>4$ since the only one step needs to be modified is showing that if $X\not\equiv 0$ then $X$ has only a finite number of positive zeros which requires the conclusions of \cite[Proposition 2]{El77} and \cite[p.273]{Sw92} and they are indeed correct for all $M>4$). Now, let $\tilde{X}_0\not\equiv 0$ be another solution of \eqref{rwevpb} with $k=0$. Then $X_0(0), \tilde{X}_0(0)\neq 0$, for otherwise $X_0$ resp. $\tilde{X}_0$ would vanish identically by \cite[Lemma 2.4]{BWW03}. So we can find $\tau\in\mathbb{R}$ such that $X_0(0)=\tau\tilde{X}_0(0)$. But then $X_0-\tau\tilde{X}_0$ also solves \eqref{rwevpb} with $k=0$ and equals zero at origin. By \cite[Lemma 2.4]{BWW03}, one has $X_0-\tau\tilde{X}_0\equiv 0$, and so $\tilde{X}_0$ is a scalar multiple of $X_0$. When $\beta=\beta_{\mathrm{FS}}(\alpha)$ which implies $q^2\lambda_1=\varpi_1$, then \eqref{rwevpb} with $k=1$ admits one solution $X_1(s)=\frac{s}{(1+s^2)^{\frac{M-2}{2}}}$, in fact, when $M$ is an integer then we can directly obtain the uniqueness of solutions (up to multiplications) by using the standard  stereographic projection as in \cite{Fe02}, furthermore, since \eqref{rwevpb} is an ODE then the uniqueness also holds for all $M>4$ with minor changes.

    {\bf We claim that when $\beta<\beta_{\mathrm{FS}}(\alpha)$ for all $k\geq 1$,  \eqref{rwevpb} does not exist nontrivial solutions, and also when $\beta=\beta_{\mathrm{FS}}(\alpha)$ for all $k\geq 2$}.

    Now, we begin to show this claim when $M$ is an integer. One easily checks the operator identity
    \[
    \left(\Delta_s-\frac{\varpi_k}{s^2}\right)(\cdot)
    =s^k\left[\frac{\partial^2}{\partial s^2}+\frac{M+2k-1}{s}\frac{\partial}{\partial s}\right](s^{-k}\cdot).
    \]
    Therefore equation \eqref{rwevpb} can be rewritten as
    \begin{align}\label{rwevpbb}
    \left(\frac{\partial^2}{\partial s^2}+\frac{M+2k-1}{s}\frac{\partial}{\partial s}\right)^2 \bar{X}_{k,i}
    & = (\tilde{2}^{**}-1)\Gamma_M(1+s^2)^{-4}\bar{X}_{k,i}
    +\left(q^2\lambda_k-\varpi_k\right)
    \bigg[\frac{2}{s^2}\bar{X}_{k,i}''
    \nonumber\\
    & \quad\quad+\frac{2(M-3)}{s^3}\bar{X}_{k,i}'
    -\frac{2(M-4)+q^2\lambda_k+\varpi_k}{s^4}\bar{X}_{k,i}\bigg].
    \end{align}
    Here we defined $\bar{X}_{k,i}\in C^\infty_0((0,\infty))$ by $\bar{X}_{k,i}(s):=s^{-k}X_{k,i}$.
    Now we consider the function $Z_{k,i}: \mathbb{R}^{M+2k}\to \mathbb{R}$ defined by $Z_{k,i}(y)=\bar{X}_{k,i}(|y|)$. So following the work of Bartsch et al. \cite{BWW03}, we deduce that
    $Z_{k,i}\in \mathcal{D}^{2,2}_0(\mathbb{R}^{M+2k})$, where $\mathcal{D}^{2,2}_0(\mathbb{R}^{M+2k})$ denotes the completion of $C^\infty_0(\mathbb{R}^{M+2k})$ with respect to the norm
    \[
    \|u\|_{\mathcal{D}^{2,2}_0(\mathbb{R}^{M+2k})}
    =\left(\int_{\mathbb{R}^{M+2k}}|\Delta u|^2 \mathrm{d}y\right)^{\frac{1}{2}},
    \]
    and $Z_{k,i}$ is a weak solution of the equation
    \begin{align}\label{rwevpbbb}
    \Delta^2 Z_{k,i}(y)& = \left(q^2\lambda_k-\varpi_k\right)
    \left(\frac{2}{s^2}Z_{k,i}''
    +\frac{2(M-3)}{s^3}Z_{k,i}'
    -\frac{2(M-4)+q^2\lambda_k+\varpi_k}{s^3}Z_{k,i}\right)
    \nonumber \\
    & \quad+(\tilde{2}^{**}-1)\Gamma_M(1+|y|^2)^{-4}Z_{k,i}(y),\quad y\in \mathbb{R}^{M+2k}.
    \end{align}
    Multiplying \eqref{rwevpbbb} by $Z_{k,i}$ and integrating in $\mathbb{R}^{M+2k}$, we have
    \begin{small}\begin{align*}%\label{rwevpbbbc}
    \|Z_{k,i}\|^2_{\mathcal{D}^{2,2}_0(\mathbb{R}^{M+2k})}
    & = (\tilde{2}^{**}-1)\Gamma_M\int_{\mathbb{R}^{M+2k}}
    (1+|y|^2)^{-4}|Z_{k,i}|^2 \mathrm{d}y
    -\left(q^2\lambda_k-\varpi_k\right)
    \Bigg\{
    2\int_{\mathbb{R}^{M+2k}}
    \frac{|\nabla Z_{k,i}|^2}{|y|^{2}}
    \mathrm{d}y
    \nonumber \\
    &\quad \quad+[(2(M-4)+q^2\lambda_k+\varpi_k)-2k(M+2k-4)]
    \int_{\mathbb{R}^{M+2k}}\frac{|Z_{k,i}|^2}
    {|y|^{4}}\mathrm{d}y
    \Bigg\}.
    \end{align*}\end{small}
    By the classical Hardy inequality,
    \begin{align*}
    \left(\frac{M+2k-4}{2}\right)^2\int_{\mathbb{R}^{M+2k}}
    \frac{|u|^2}{|y|^{4}}\mathrm{d}y
    \leq \int_{\mathbb{R}^{M+2k}}
    \frac{|\nabla u|^2}{|y|^{2}}\mathrm{d}y, \quad \mbox{for all}\quad u\in C^\infty_0(\mathbb{R}^{M+2k}),
    \end{align*}
    %and the inequality is strict for nontrivial functions.
    we deduce that
    \begin{align}\label{rwevpbbbcb}
    \|Z_{k,i}\|^2_{\mathcal{D}^{2,2}_0(\mathbb{R}^{M+2k})}
    & \leq (\tilde{2}^{**}-1)\Gamma_M\int_{\mathbb{R}^{M+2k}}
    (1+|y|^2)^{-4}|Z_{k,i}|^2 \mathrm{d}y
    \nonumber\\
    & \quad-\left(q^2\lambda_k-\varpi_k\right)\xi_k
    \int_{\mathbb{R}^{M+2k}}\frac{|Z_{k,i}|^2}
    {|y|^{4}}\mathrm{d}y.
    \end{align}
    Here
    \[
    \xi_k:=\frac{(M+2k-4)^2}{2}+(2(M-4)+q^2\lambda_k+\varpi_k)
    -2k(M+2k-4).
    \]
    From \cite[(2.10)]{BWW03},
    \begin{align}\label{rwevpbbbcbi}
    \|u\|^2_{\mathcal{D}^{2,2}_0(\mathbb{R}^{M+2k})}
    \geq & \Gamma_{M+2k}\int_{\mathbb{R}^{M+2k}}
    \frac{|u(y)|^2}{(1+|y|^2)^{4}} \mathrm{d}y,\quad \mbox{for all}\quad u\in \mathcal{D}^{2,2}_0(\mathbb{R}^{M+2k}),
    \end{align}
    then combining with \eqref{rwevpbbbcb} and \eqref{rwevpbbbcbi} we deduce
    \begin{align*}
    & \left[(\tilde{2}^{**}-1)\Gamma_M-\Gamma_{M+2k}\right]
    \int_{\mathbb{R}^{M+2k}}
    (1+|y|^2)^{-4}|Z_{k,i}|^2 \mathrm{d}y
    \geq \left(q^2\lambda_k-\varpi_k\right)\xi_k
    \int_{\mathbb{R}^{M+2k}}\frac{|Z_{k,i}|^2}
    {|y|^{4}}\mathrm{d}y.
    \end{align*}
    Since
    \[
    (\tilde{2}^{**}-1)\Gamma_M\leq \Gamma_{M+2k},\quad \mbox{for all}\quad k\geq 1,
    \]
    and
    \[
    (\tilde{2}^{**}-1)\Gamma_M< \Gamma_{M+2k},\quad \mbox{for all}\quad k\geq 2,
    \]
    where $\Gamma_M$ is defined in \eqref{defgn}, then it holds that
    \begin{align}\label{rwevpbbbcbcf}
    \left(q^2\lambda_k-\varpi_k\right)\xi_k
    \int_{\mathbb{R}^{M+2k}}\frac{|Z_{k,i}|^2}
    {|y|^{4}}\mathrm{d}y\leq 0.
    \end{align}
    Note that
    \begin{align*}
    \xi_k
    & \geq \frac{(M+2k-4)^2}{2}+(2(M-4)+2\varpi_k)
    -2k(M+2k-4) \\
    & = \frac{M^2}{2}+2(M-2)(k-2)+2(M-4)>0,
    \end{align*}
    for all $k\geq 1$, then we conclude from \eqref{rwevpbbbcbcf} that $Z_{k,i}\equiv 0$ for all $k\geq 1$ if $\beta<\beta_{\mathrm{FS}}(\alpha)$, and $Z_{k,i}\equiv 0$ for all $k\geq 2$ if $\beta=\beta_{\mathrm{FS}}(\alpha)$.
    Since \eqref{rwevpb} is ODE, even if $M$ is not an integer we readily see that the above conclusion remains true. Our claim is proved.

    To sum up, let us turn back to (\ref{p2c}) we obtain the solutions that, if $\beta=\beta_{\mathrm{FS}}(\alpha)$ then \eqref{p2c} only admits
    \begin{equation*}%\label{pye}
    \phi_0(r)=\frac{1-r^{2+\beta-\alpha}}
    {(1+r^{2+\beta-\alpha})^{\frac{N-2+\alpha}{2+\beta-\alpha}}},\quad
    \phi_1(r)=\frac{r^{\frac{2+\beta-\alpha}{2}}}
    {(1+r^{2+\beta-\alpha})^{\frac{N-2+\alpha}{2+\beta-\alpha}}},
    \end{equation*}
    otherwise, \eqref{p2c} only admits $\phi_0$.
    That is, if $\beta=\beta_{\mathrm{FS}}(\alpha)$, the space of solutions of (\ref{Pwhlp}) has dimension $(1+N)$ and is spanned by
    \begin{equation*}
    Z_{0}(x)=\frac{1-|x|^{2+\beta-\alpha}}{(1+|x|^{2+\beta-\alpha})
    ^\frac{N-2+\alpha}{2+\beta-\alpha}},\quad Z_{i}(x)=\frac{|x|^{\frac{2+\beta-\alpha}{2}}}{(1+|x|^{2+\beta-\alpha})
    ^\frac{N-2+\alpha}{2+\beta-\alpha}}\cdot \frac{x_i}{|x|},\quad i=1,\ldots,N.
    \end{equation*}
    Otherwise the space of solutions of (\ref{Pwhlp}) has dimension one and is spanned by $Z_0$, and note that $Z_0\thicksim \frac{\partial U_{\lambda}}{\partial \lambda}|_{\lambda=1}$ in this case we say $U$ is non-degenerate. The proof of Theorem \ref{thmpwhl} is now completed.
    \qed

\vskip0.25cm

Now, based the proof of Theorem \ref{thmpwhl}, we are ready to prove the symmetry breaking phenomenon by taking the same arguments as those in \cite{FS03} of Felli and Schneider. In order to shorten formulas, for each $u\in \mathcal{D}^{2,2}_{\alpha,\beta}(\mathbb{R}^N)$, we denote
    \begin{equation}\label{def:norm}
    \|u\|: =\left(\int_{\mathbb{R}^N}|x|^{-\beta}|\mathrm{div} (|x|^{\alpha}\nabla u)|^2 \mathrm{d}x\right)^{\frac{1}{2}},
    \quad \|u\|_*: =\left(\int_{\mathbb{R}^N}|x|^\beta|u|^{p^*_{\alpha,\beta}} \mathrm{d}x\right)^{\frac{1}{p^*_{\alpha,\beta}}}.
    \end{equation}

\subsection{Proof of Theorem \ref{thmmr}}\label{subsectmrt}
    We follow the arguments in the proof of \cite[Corollary 1.2]{FS03}. We define the functional $\mathcal{I}$ on $\mathcal{D}^{2,2}_{\alpha,\beta}(\mathbb{R}^N)$ by the right hand side of \eqref{ckns}, i.e.,
    \begin{align}\label{defFu}
    \mathcal{I}(u):=\frac{\|u\|^2}{\|u\|_*^2},\quad u\in \mathcal{D}^{2,2}_{\alpha,\beta}(\mathbb{R}^N)\setminus\{0\}.
    \end{align}
    Define also the energy functional of equation \eqref{Pwh} as
    \[
    \mathcal{J}(u):=\frac{1}{2}\|u\|^2
    -\frac{1}{p^*_{\alpha,\beta}}\|u\|^{p^*_{\alpha,\beta}}_*, \quad u\in \mathcal{D}^{2,2}_{\alpha,\beta}(\mathbb{R}^N).
    \]
    From Theorem \ref{thmpwh} we know $U$ given in \eqref{defula} is a critical of $\mathcal{J}$, thus $\langle \mathcal{J}'(U),\varphi\rangle=0$ for all $\varphi\in \mathcal{D}^{2,2}_{\alpha,\beta}(\mathbb{R}^N)$ and $\|U\|^2=\|U\|_*^{p^*_{\alpha,\beta}}$.
    Then we obtain for all $\varphi_1,\varphi_2\in \mathcal{D}^{2,2}_{\alpha,\beta}(\mathbb{R}^N)$, it holds that
    \begin{align*}
    \langle \mathcal{I}'(U),\varphi_1\rangle
    = \frac{2}{\|U\|_*^2}
    \langle \mathcal{J}'(U),\varphi_1\rangle
    =0,
    \end{align*}
    and
    \begin{align*}
    \langle \mathcal{I}''(U)\varphi_1,\varphi_2\rangle
    & = \frac{2}{\|U\|_*^2}\langle \mathcal{J}''(U)\varphi_1,\varphi_2\rangle
    + \frac{2(p^*_{\alpha,\beta}-2)}{\|U\|_*^{p^*_{\alpha,\beta}+2}}
    \int_{\mathbb{R}^N}|x|^\beta U^{p^*_{\alpha,\beta}-1} \varphi_1 \mathrm{d}x
    \int_{\mathbb{R}^N}|x|^\beta U^{p^*_{\alpha,\beta}-1} \varphi_2 \mathrm{d}x
    \\
    & \quad-\frac{4}{\|U\|_*^{p^*_{\alpha,\beta}+2}}
    \int_{\mathbb{R}^N}|x|^{-\beta}
    \mathrm{div}(|x|^{\alpha}\nabla U) \mathrm{div}(|x|^{\alpha}\nabla \varphi_1)\mathrm{d}x
    \int_{\mathbb{R}^N}|x|^\beta
    U^{p^*_{\alpha,\beta}-1}\varphi_2 \mathrm{d}x
    \\
    & \quad-\frac{4}{\|U\|_*^{p^*_{\alpha,\beta}+2}}
    \int_{\mathbb{R}^N}|x|^{-\beta}
    \mathrm{div}(|x|^{\alpha}\nabla U) \mathrm{div}(|x|^{\alpha}\nabla \varphi_2)\mathrm{d}x
    \int_{\mathbb{R}^N}|x|^\beta
    U^{p^*_{\alpha,\beta}-1}\varphi_1 \mathrm{d}x.
    \end{align*}
    From the proof of Theorem \ref{thmpwhl}, we know that $X_1(s)=\frac{s}{(1+|s|^2)^{\frac{M-2}{2}}}$ with $M=\frac{2(N+\beta)}{2+\beta-\alpha}>4$ is a solution of
    \begin{eqnarray*}
    \left\{ \arraycolsep=1.5pt
       \begin{array}{ll}
        X_1''+\frac{M-1}{s}X_1'-\frac{M-1}{s^2}X_1+Y=0 \quad \mbox{in}\quad s\in(0,\infty),\\[3mm]
        Y''+\frac{M-1}{s}Y'-\frac{M-1}{s^2}Y
        +\frac{(M+4)(M-2)M(M+2)}{(1+s^2)^4}X_1=0 \quad \mbox{in}\quad s\in(0,\infty).
        \end{array}
    \right.
    \end{eqnarray*}
    Therefore, when $\alpha>0$ and $\beta_{\mathrm{FS}}(\alpha)<\beta<\frac{N}{N-2}\alpha$ which imply $\mu:=\left(\frac{2}{2+\beta-\alpha}\right)^2(N-1)<M-1$, then we deduce that
    \begin{align*}
    \langle\mathcal{J}''(U)Z_{i},Z_{i}\rangle
    & = \frac{\omega_{N-1}}{N}\left(\frac{2+\beta-\alpha}{2}\right)^3
    [\mu-(M-1)]
    \bigg\{
    2\int^\infty_0(X_1'(s))^2 s^{M-4}\mathrm{d}s
    \\
    & \quad \quad+[(2M-5)+\mu]\int^\infty_0(X_1(s))^2 s^{M-5}\mathrm{d}s
    \bigg\}
    < 0,
    \end{align*}
    where
    $Z_{i}(x)=X_1(|x|^{\frac{2+\beta-\alpha}{2}})\frac{x_i}{|x|}$ for some $i\in \{1,2,\ldots,N\}$ same as in \eqref{defaezki}. Then we infer
    \[
    \langle \mathcal{I}''(U)Z_{i},Z_{i}\rangle<0,
    \]
    due to $\int_{\mathbb{R}^N}|x|^\beta U^{p^*_{\alpha,\beta}-1} Z_{i}\mathrm{d}x=0$.
    Consequently, $\mathcal{S}$ is strictly small than $\mathcal{I}(U)=\mathcal{S}_r$. Then no extremal functions of $\mathcal{S}$ are symmetry radial, the proof of Theorem \ref{thmmr} is now completed.
    \qed

\section{{\bfseries A symmetry result}}\label{sectps}

In this section, we will show a symmetry result for the seconde-order (CKN) type inequality \eqref{ckn2n} when $p^{*}_{\alpha,\beta}=2^{**}$ and $2-N<\alpha<0$, and give the proof of Theorem \ref{thm2ps}. This bases on a crucial Rellich-Sobolev type inequality with explicit form of extremal functions which was established by Dan et al. \cite{DMY20} given as in \eqref{RSi}.

\subsection{Proof of Theorem \ref{thm2ps}}\label{subsectsr}
    For each $u\in \mathcal{D}^{2,2}_{\alpha,\frac{N}{N-2}\alpha}(\mathbb{R}^N)$, let us make the change
    \[
    u(x)=|x|^{\eta}v(x) \quad\mbox{with}\quad
    \eta=-\frac{N-4}{2(N-2)}\alpha.
    \]
    A direct calculation indicates
    \begin{align}\label{psle}
    \int_{\mathbb{R}^N}
    |x|^{\frac{N}{N-2}\alpha}|u|^{2^{**}} \mathrm{d}x
    =\int_{\mathbb{R}^N}
    |v|^{2^{**}} \mathrm{d}x,
    \end{align}
    and
    \begin{align}\label{psny}
    & \int_{\mathbb{R}^N}|x|^{-\frac{N}{N-2}\alpha}|\mathrm{div} (|x|^{\alpha}\nabla u)|^2 \mathrm{d}x
    \nonumber\\
    & =\int_{\mathbb{R}^N}|x|^{-\frac{N}{N-2}\alpha}
    [|x|^{\eta+\alpha}\Delta v +(2\eta+\alpha)|x|^{\alpha+\eta-2}(x\cdot \nabla v)
    +\eta(N+\alpha+\eta-2)|x|^{\alpha+\eta-2}v]^2 \mathrm{d}x
    \nonumber\\
    & = \int_{\mathbb{R}^N}
    [\Delta v +(2\eta+\alpha)|x|^{-2}(x\cdot \nabla v)
    +\eta(N+\alpha+\eta-2)|x|^{-2}v]^2 \mathrm{d}x
    \nonumber\\
    & = \int_{\mathbb{R}^N}|\Delta v|^2 \mathrm{d}x
    +\eta^2(N+\alpha+\eta-2)^2\int_{\mathbb{R}^N}|x|^{-4}v^2\mathrm{d}x
    +(2\eta+\alpha)^2\int_{\mathbb{R}^N}|x|^{-4}(x\cdot \nabla v)^2\mathrm{d}x
    \nonumber\\
    &\quad +2(2\eta+\alpha)\int_{\mathbb{R}^N}|x|^{-2}(x\cdot \nabla v)\Delta v\mathrm{d}x
    +2\eta(N+\alpha+\eta-2)\int_{\mathbb{R}^N}|x|^{-2}v\Delta v\mathrm{d}x
    \nonumber\\
    &\quad +2(2\eta+\alpha)\eta(N+\alpha+\eta-2)
    \int_{\mathbb{R}^N}|x|^{-4}v(x\cdot \nabla v)\mathrm{d}x,
    \end{align}
    due to $2(\alpha+\eta)+4-\frac{N-4}{N-2}\alpha=0$. In order to calculate those integrals in last term, we need the following identity:
    \begin{align}\label{eqi}
    (N-4)\int_{\mathbb{R}^N}
    |x|^{-2}|\nabla v|^2\mathrm{d}x
    =
    2\int_{\mathbb{R}^N}
    (x\cdot\nabla v)\mathrm{div}(|x|^{-2}\nabla v)
    \mathrm{d}x.%,\quad \mbox{for all}\quad v\in \mathcal{D}^{2,2}_0(\mathbb{R}^N).
    \end{align}
    In order to show the claim \eqref{eqi}, let $v(x)=|x|w(x)$ then
    \begin{align*}
    \int_{\mathbb{R}^N}
    |x|^{-2}|\nabla v|^2\mathrm{d}x
    =
    \int_{\mathbb{R}^N}
    [|\nabla w|^2
    +|x|^{-2}w^2+
    2|x|^{-2}w(x\cdot\nabla w)]
    \mathrm{d}x,
    \end{align*}
    and
    \begin{align*}
    \int_{\mathbb{R}^N}
    (x\cdot\nabla v)\mathrm{div}(|x|^{-2}\nabla v)
    \mathrm{d}x
    =
    \int_{\mathbb{R}^N}
    \left[
    (N-3)|x|^{-2}(w^2+w(x\cdot\nabla w))
    + \Delta w(w+ x\cdot\nabla w)
    \right]
    \mathrm{d}x.
    \end{align*}
    Therefore, it is easy to verify that \eqref{eqi} is equivalent to
    \begin{align}\label{eqib}
    (N-2)\int_{\mathbb{R}^N}
    |\nabla w|^2\mathrm{d}x
    =
    2\int_{\mathbb{R}^N}
    \Delta w(x\cdot\nabla w)
    \mathrm{d}x.
    \end{align}
    Note that
    \begin{align*}
    \int_{\mathbb{R}^N}
    \Delta w(x\cdot\nabla w)
    \mathrm{d}x
    & = -\int_{\mathbb{R}^N}
    \nabla w\cdot \nabla(x\cdot\nabla w)
    \mathrm{d}x
    \\
    & = -\int_{\mathbb{R}^N}
    \sum^{N}_{i=1}\frac{\partial w}{\partial x_i}
    \left(\frac{\partial w}{\partial x_i}
    +\sum^{N}_{j=1}x_j\frac{\partial^2 w}{\partial x_j\partial x_i}\right)
    \mathrm{d}x
    \\
    & = -\int_{\mathbb{R}^N}|\nabla w|^2\mathrm{d}x
    -\int_{\mathbb{R}^N}
    \sum^{N}_{i=1}\sum^{N}_{j=1}\frac{\partial w}{\partial x_i}x_j\frac{\partial^2 w}{\partial x_j\partial x_i}
    \mathrm{d}x.
    \end{align*}
    By using the partial integration method, we have
    \begin{align*}
    \int_{\mathbb{R}^N}
    \sum^{N}_{i=1}\sum^{N}_{j=1}\frac{\partial w}{\partial x_i}x_j\frac{\partial^2 w}{\partial x_j\partial x_i}
    \mathrm{d}x
    & = \int_{\mathbb{R}^{N-1}}\int_{\mathbb{R}}
    \sum^{N}_{i=1}\sum^{N}_{j=1}\frac{\partial w}{\partial x_i}x_j
    \mathrm{d}\frac{\partial w}{\partial x_i}\mathrm{d}x'
    = -\frac{N}{2}\int_{\mathbb{R}^N}|\nabla w|^2\mathrm{d}x,
    \end{align*}
    where $x=(x',x_j)\in \mathbb{R}^{N-1}\times\mathbb{R}$, thus \eqref{eqib} holds then \eqref{eqi} also holds. From \eqref{eqi}, we obtain
    \begin{align}\label{psd}
    (N-4)\int_{\mathbb{R}^N}
    |x|^{-2}|\nabla v|^2\mathrm{d}x
    = 2\int_{\mathbb{R}^N}
    |x|^{-2}(x\cdot\nabla v)\Delta v
    \mathrm{d}x-4\int_{\mathbb{R}^N}|x|^{-4}(x\cdot\nabla v)^2
    \mathrm{d}x,
    \end{align}
    and by Green's formula
    \begin{align}\label{psp}
    \int_{\mathbb{R}^N}|x|^{-2}v\Delta v\mathrm{d}x
    & =-\int_{\mathbb{R}^N}\nabla(|x|^{-2}v)\cdot\nabla v\mathrm{d}x
    \nonumber\\
    & = -\int_{\mathbb{R}^N}|x|^{-2}|\nabla v|^2\mathrm{d}x
    +2\int_{\mathbb{R}^N}|x|^{-4}v(x\cdot\nabla v)
    \mathrm{d}x,
    \end{align}
    and by Divergence formula
    \begin{align}\label{psdf}
    \int_{\mathbb{R}^N}|x|^{-4}v(x\cdot \nabla v)\mathrm{d}x
    & =-\int_{\mathbb{R}^N}v\mathrm{div}(|x|^{-4}xv)\mathrm{d}x
    \nonumber\\
    & = -(N-4)\int_{\mathbb{R}^N}|x|^{-4}v^2\mathrm{d}x
    -\int_{\mathbb{R}^N}|x|^{-4}v(x\cdot \nabla v)\mathrm{d}x
    \nonumber\\
    & = -\frac{N-4}{2}\int_{\mathbb{R}^N}|x|^{-4}v^2\mathrm{d}x.
    \end{align}
    Therefore, from \eqref{psd}, \eqref{psp} and \eqref{psdf} we deduce
    \begin{align*}
    & (2\eta+\alpha)^2\int_{\mathbb{R}^N}|x|^{-4}(x\cdot \nabla v)^2\mathrm{d}x+2(2\eta+\alpha)\int_{\mathbb{R}^N}|x|^{-2}(x\cdot \nabla v)\Delta v\mathrm{d}x
    \\
    &\quad+2\eta(N+\alpha+\eta-2)\int_{\mathbb{R}^N}|x|^{-2}v\Delta v\mathrm{d}x
    +2(2\eta+\alpha)\eta(N+\alpha+\eta-2)
    \int_{\mathbb{R}^N}|x|^{-4}v(x\cdot \nabla v)\mathrm{d}x
    \\
    &= -(N-4)\eta(N+\alpha+\eta-2)(2\eta+\alpha+2)
    \int_{\mathbb{R}^N}|x|^{-4}v^2\mathrm{d}x
    \\
    &\quad+ (2\eta+\alpha)(2\eta+\alpha+4)\int_{\mathbb{R}^N}|x|^{-4}(x\cdot\nabla v)^2 \mathrm{d}x
    \\
    &\quad +[(2\eta+\alpha)(N-4)-2\eta(N+\alpha+\eta-2)]
    \int_{\mathbb{R}^N}|x|^{-2}|\nabla v|^2\mathrm{d}x
    \\
    & \geq -(N-4)\eta(N+\alpha+\eta-2)(2\eta+\alpha+2)
    \int_{\mathbb{R}^N}|x|^{-4}v^2\mathrm{d}x
    \\
    &\quad +[(2\eta+\alpha)(N+2\eta+\alpha)-2\eta(N+\alpha+\eta-2)]
    \int_{\mathbb{R}^N}|x|^{-2}|\nabla v|^2\mathrm{d}x,
    \end{align*}
    due to $(x\cdot\nabla v)^2\leq |x|^2|\nabla v|^2$ and $(2\eta+\alpha)(2\eta+\alpha+4)=\frac{4\alpha}{N-2}
    \left(2+\frac{\alpha}{N-2}\right)<0$ (note that the assumption $\alpha<0$ plays a crucial role), furthermore, the equality holds if and only if $v$ is radial.
    Thus from \eqref{psny} we have
    \begin{align}\label{psnb}
    & \int_{\mathbb{R}^N}|x|^{-\frac{N}{N-2}\alpha}|\mathrm{div} (|x|^{\alpha}\nabla u)|^2 \mathrm{d}x
    \geq \int_{\mathbb{R}^N}|\Delta v|^2 \mathrm{d}x
    \nonumber\\
    &\quad+ [(2\eta+\alpha)(N+2\eta+\alpha)-2\eta(N+\alpha+\eta-2)]
    \int_{\mathbb{R}^N}|x|^{-2}|\nabla v|^2\mathrm{d}x
    \nonumber\\
    &\quad+ \left[\eta^2(N+\alpha+\eta-2)^2
    -(N-4)\eta(N+\alpha+\eta-2)(2\eta+\alpha+2)\right]
    \int_{\mathbb{R}^N}|x|^{-4}v^2\mathrm{d}x,
    \end{align}
    and the equality holds if and only if $v$ is radial (so does $u$ due to $u(x)=|x|^\eta v(x)$). Now, let us make the change
    \[
    \alpha=\frac{2-N}{N-4}\mu,
    \]
    which implies $0<\mu<N-4$ due to $2-N<\alpha<0$, then from $\eta=-\frac{N-4}{2(N-2)}\alpha$ we have
    \begin{align*}
    (2\eta+\alpha)(N+2\eta+\alpha)-2\eta(N+\alpha+\eta-2)
    & =\frac{\alpha[2(N-2)+\alpha]}{2(N-2)^2}(N^2-4N+8)
    \\
    & = -\frac{N^2-4N+8}{2(N-4)^2}\mu[2(N-4)-\mu]
    =:-C_{\mu,1},
    \end{align*}
    and
    \begin{align*}
    & \eta^2(N+\alpha+\eta-2)^2-(N-4)\eta(N+\alpha+\eta-2)(2\eta+\alpha+2)
    \\
    &= \frac{N^2}{16(N-4)^2}\mu^2[2(N-4)-\mu]^2
    -\frac{N-2}{2}\mu[2(N-4)-\mu]
    =:C_{\mu,2},
    \end{align*}
    Therefore, \eqref{psnb} indicates
    \begin{align*}
    \int_{\mathbb{R}^N}|x|^{-\frac{N}{N-2}\alpha}|\mathrm{div} (|x|^{\alpha}\nabla u)|^2 \mathrm{d}x
    & \geq \int_{\mathbb{R}^N}|\Delta v|^2 \mathrm{d}x
    - C_{\mu,1}\int_{\mathbb{R}^N}|x|^{-2}|\nabla v|^2\mathrm{d}x
    + C_{\mu,2}\int_{\mathbb{R}^N}|x|^{-4}v^2\mathrm{d}x,
    \end{align*}
    and equality holds only if $v$ is radial. Therefore from the Rellich-Sobolev type inequality \eqref{RSi} established by Dan et al. in \cite[Theorem 1.6]{DMY20}, and \eqref{psle} we deduce
    \begin{align*}
    \int_{\mathbb{R}^N}|x|^{-\frac{N}{N-2}\alpha}|\mathrm{div} (|x|^{\alpha}\nabla u)|^2 \mathrm{d}x
    \geq \left(1+\frac{\alpha}{N-2}\right)^{4-\frac{4}{N}}\mathcal{S}_0
    \left(\int_{\mathbb{R}^N}
    |x|^{\frac{N}{N-2}\alpha}|u|^{2^{**}} \mathrm{d}x\right)^{\frac{2}{2^{**}}},
    \end{align*}
    and equality holds if and only if
    \[
    u(x)=c|x|^\eta|x|^{-\frac{\mu}{2}}
    \left(\nu+|x|^{2(1-\frac{\mu}{N-4})}
    \right)^{-\frac{N-4}{2}}
    =A
    (\nu+|x|^{2+\frac{2\alpha}{N-2}})
    ^{-\frac{N-4}{2}},\quad \mbox{for}\ A\in\mathbb{R},\ \nu>0,
    \]
    that is, $u(x)=A\lambda^{\frac{N-4}{2}(1+\frac{\alpha}{N-2})}U(\lambda x)$ for $A\in\mathbb{R}$ and $\lambda>0$, where
    \begin{align*}
    U(x)=C_{N,\alpha,\frac{N}{N-2}\alpha}
    (1+|x|^{2+\frac{2\alpha}{N-2}})
    ^{-\frac{N-4}{2}}
    \end{align*}
    is given as in Theorem \ref{thmpwh} replacing $\beta$ by $\frac{N}{N-2}\alpha$.
    Now the proof of Theorem \ref{thm2ps} is completed.
    \qed

\medskip

\noindent{\bfseries Statements and Declarations}
The authors declare that they have no conflict of interest.

\medskip

\noindent{\bfseries Data availability}
\noindent No data was used for the research described in the article.

\medskip

\noindent{\bfseries Acknowledgements}
The research has been supported by National Natural Science Foundation of China (No. 12371121).

    \end{document}